\documentclass{article}
\usepackage{stmaryrd}
\usepackage{mathrsfs}
\usepackage[centertags]{amsmath}
\usepackage{amsfonts, dsfont}
\usepackage{amssymb}
\usepackage{amsthm}
\usepackage{newlfont}
\usepackage{bezier,  amsxtra, amsbsy, amsgen, amsopn, amstext}

\input xy
\xyoption{all}

\theoremstyle{plain}
\newtheorem{thm}{Theorem}[section]
\newtheorem{cor}[thm]{Corollary}
\newtheorem{lem}[thm]{Lemma}
\newtheorem{prop}[thm]{Proposition}

\theoremstyle{definition}
\newtheorem{defn}[thm]{Definition}
\newtheorem{remark}[thm]{Remark}
\newtheorem*{ack}{Acknowledgments}

\newcommand{\bd}{\begin{defn}}
\newcommand{\ed}{\end{defn}}
\newcommand{\bl}{\begin{lem}}
\newcommand{\el}{\end{lem}}
\newcommand{\bp}{\begin{prop}}
\newcommand{\ep}{\end{prop}}
\newcommand{\bt}{\begin{thm}}
\newcommand{\et}{\end{thm}}
\newcommand{\bc}{\begin{cor}}
\newcommand{\ec}{\end{cor}}
\newcommand{\br}{\begin{remark}}
\newcommand{\er}{\end{remark}}
\newcommand{\bdi}{\begin{diagram}}
\newcommand{\edi}{\end{diagram}}
\newcommand{\beq}{\begin{eqn}}
\newcommand{\eeq}{\end{eqn}}
\newcommand{\ba}{\begin{array}}
\newcommand{\ea}{\end{array}}
\newcommand{\bpf}{\begin{proof}}
\newcommand{\epf}{\end{proof}}

\newcommand{\Q}{\mathds{Q}}
\newcommand{\Zp}{\mathds{Z}_{p}}
\newcommand{\Qp}{\mathds{Q}_{p}}
\newcommand{\al}{\alpha}
\newcommand{\be}{\beta}

\newcommand{\Si}{\Sigma}
\newcommand{\Ga}{\Gamma}
\newcommand{\ga}{\gamma}
\newcommand{\La}{\Lambda}
\newcommand{\la}{\lambda}

\newcommand{\Op}{\mathcal{O}}

\newcommand{\M}{\mathfrak{M}}

\DeclareMathOperator{\Gal}{Gal} \DeclareMathOperator{\Hom}{Hom}
\DeclareMathOperator{\Ext}{Ext} 
 \DeclareMathOperator{\rank}{rank}
\newcommand{\ot}{\otimes}

\newcommand{\ilim}{\displaystyle \mathop{\varinjlim}\limits}
\newcommand{\plim}{\displaystyle \mathop{\varprojlim}\limits}
\newcommand{\im}{\mathrm{im}\,}
\newcommand{\coker}{\mathrm{coker}\,}

\newcommand{\cyc}{\mathrm{cyc}}

\newcommand{\lra}{\longrightarrow}
\newcommand{\tha}{\twoheadrightarrow}

\newcommand{\sbs}{\subseteq}

\newcommand{\ps}[1]{\llbracket #1 \rrbracket}

\begin{document}

\title{A remark on the $\mathfrak{M}_H(G)$-conjecture and Akashi series}
\author{Meng Fai Lim\footnote{Department of Mathematics, University of Toronto, 40 St. George St.,
Toronto, Ontario, Canada M5S 2E4}}
\date{}
\maketitle

\begin{abstract} \footnotesize
\noindent  In this article, we give a criterion for the dual Selmer
group of an elliptic curve which has either good ordinary reduction
or multiplicative reduction at every prime above $p$ to satisfy the
$\M_H(G)$-conjecture. As a by-product of our calculations, we are
able to define the Akashi series of the dual Selmer groups assuming
the conjectures of Mazur and Schneider. Previously, the Akashi
series are defined under the stronger assumption that the dual
Selmer group satisfies the $\M_H(G)$-conjecture. We then establish a
criterion for the vanishing of the dual Selmer groups using the
Akashi series. We will apply this criterion to prove some results on
the characteristic elements of the dual Selmer groups. Our methods
in this paper are inspired by the work of Coates-Schneider-Sujatha
and can be extended to the Greenberg Selmer groups attached to other
ordinary representations, for instance, those coming from an
$p$-ordinary modular form.



\end{abstract}

\section{Introduction}

Throughout the paper, $p$ will always denote an odd prime. Let $G$
be a compact $p$-adic Lie group with a closed normal subgroup $H$
such that $G/H\cong \Zp$. We denote $\mathfrak{M}_H(G)$ to be the
category of a special class of $\Zp\ps{G}$-torsion modules, namely
those finitely generated left $\Zp\ps{G}$-modules $M$ such that
$M/M(p)$ is finitely generated over $\Zp\ps{H}$; here $M(p)$ is the
$\Zp\ps{G}$-submodule of $M$ consisting of elements of $M$ which are
annihilated by some power of $p$. In this paper, we give a criterion
for the dual Selmer group of an elliptic curve to be in
$\mathfrak{M}_H(G)$. Such a criterion is already known in certain
cases when $G$ is the Galois group of a $p$-adic Lie extension of
dimension 2 and when $G$ is the Galois group of the $p$-adic Lie
extension obtained by adjoining all the $p$-division points of the
elliptic curve. The criterion we obtain here can therefore be viewed
as a generalization of those. We remark that the knowledge of such
dual Selmer groups lying in $\mathfrak{M}_H(G)$ is crucial in the
formulation of the main conjectures of non-commutative Iwasawa
theory, as one requires this property in order to be able to define
suitable characteristic elements for these modules in certain
relative $K_0$-groups (see \cite{CFKSV}).

It will follow from our calculations that the Selmer group has
better properties than a general module in $\M_H(G)$ under certain
assumptions (see Proposition \ref{Selmer fg}). This in turn allows
us to define the ``Akashi series" of the dual Selmer group.
Previously, these are defined under the assumption that the dual
Selmer group lies in $\mathfrak{M}_H(G)$ (see \cite[Section
3]{CFKSV} and \cite[Section 2]{Z11}). We then establish a criterion
for the vanishing of the dual Selmer group of an elliptic curve with
good ordinary reduction at all primes above $p$ under certain mild
assumptions (see Theorems \ref{Selmer vanishing} and \ref{twist
Selmer vanishing}). This is then applied to establish one of the
cases of \cite[Conjecture 4.8]{CFKSV} for this dual Selmer group
(see Propositions \ref{conjecture 4.8} and \ref{conjecture 4.8(a)}).

Our methods in this paper follow the approach in \cite[Section
2]{CSS} with supplementary ideas and techniques from \cite{Ho}. In
the final section of the paper, we will mention briefly how these
techniques can be extended to consider Greenberg's Selmer groups of
more general ordinary representations.

\section{Algebraic Preliminaries} \label{algebra}

In this section, we recall some algebraic preliminaries that will be
required in the later part of the paper. The materials presented
here base mainly on \cite{Ho}. Let $G$ be a compact pro-$p$ $p$-adic
Lie group without $p$-torsion. It is well known that $\Zp\llbracket
G\rrbracket$ is an Auslander regular ring (cf. \cite[Theorems
3.26]{V02}). Furthermore, the ring $\Zp\ps{G}$ has no zero divisors
(cf.\ \cite{Neu}), and therefore, admits a skew field $K(G)$ which
is flat over $\Zp\ps{G}$ (see \cite[Chapters 6 and 10]{GW} or
\cite[Chapter 4, \S 9 and \S 10]{Lam}). If $M$ is a finitely
generated $\Zp\ps{G}$-module, we define the $\Zp\ps{G}$-rank of $M$
to be
$$ \rank_{\Zp\ps{G}}M  = \dim_{K(G)} K(G)\ot_{\Zp\ps{G}}M. $$
 We say that the module $M$ is a
\textit{torsion} $\Zp\ps{G}$-module if $\rank_{\Zp\ps{G}} M = 0$.
Now suppose that $N$ is a $\mathbb{F}_p\ps{G}$-module. We then
define its $\mathbb{F}_p\ps{G}$-rank by
 $$ \rank_{\mathbb{F}_p\ps{G}} N=
 \displaystyle\frac{\rank_{\mathbb{F}_p\ps{G_0}}N}{|G:G_0|}, $$
where $G_0$ is an open normal uniform pro-$p$ subgroup of $G$. This
is integral and independent of the choice of $G_0$ (see
\cite[Proposition 1.6]{Ho}). Similarly, we will say that that $N$ is
a \textit{torsion} $\mathbb{F}_p\ps{G}$-module if
$\rank_{\mathbb{F}_p\ps{G}}N = 0$.

For a general compact $p$-adic Lie group $G$ and a finitely
generated $\Zp\ps{G}$-module $M$, we say that $M$ is a torsion
$\Zp\ps{G}$-module if there exists an open normal uniform pro-$p$
subgroup $G_0$ of $G$ such that $M$ is a torsion
$\Zp\ps{G_0}$-module in the above sense. We will also make use of a
well-known equivalent definition for $M$ to be torsion
$\Zp\ps{G}$-module, namely: $\Hom_{\Zp\ps{G}}(M,\Zp\ps{G})=0$.  The
notion of a torsion $\mathbb{F}_p\ps{G}$-module for a general
compact $p$-adic Lie group $G$ is extended in a similar fashion.

For a given finitely generated $\Zp\ps{G}$-module $M$, we denote
$M(p)$ to be the $\Zp\ps{G}$-submodule of $M$ consisting of elements
of $M$ which are annihilated by some power of $p$. Since the ring
$\Zp\ps{G}$ is Noetherian, the module $M(p)$ is finitely generated
over $\Zp\ps{G}$. Therefore, one can find an integer $r\geq 0$ such
that $p^r$ annihilates $M(p)$. Following \cite[Formula (33)]{Ho}, we
define
  \[\mu_G(M) = \sum_{i\geq 0}\rank_{\mathbb{F}_p\ps{G}}\big(p^i
   M(p)/p^{i+1}\big). \]
(For another alternative, but equivalent, definition, see
\cite[Definition 3.32]{V02}.) By the above discussion, the sum on
the right is a finite one. Also, it is clear from the definition
that $\mu_G(M) = \mu_G(M(p))$. Finally, it is not difficult to see
that this definition coincides with the classical notion of the
$\mu$-invariant for $\Ga$-modules when $G=\Ga$. We now record
certain properties of this invariant.

\bl \label{mu lemma} Let $G$ be a compact pro-$p$ $p$-adic Lie group
without $p$-torsion. Then we have the following statements.
\begin{enumerate}
\item[$(a)$] For every finitely generated $\Zp\ps{G}$-module $M$, one has
$$\mu_G(M) = \sum_{i\geq 0} (-1)^i\mathrm{ord}_p\big(H_i(G,M(p))\big).$$

\item[$(b)$]  Suppose that $G$ has a closed normal subgroup $H$ such that
$G/H\cong \Zp$. If $M$ is a $\Zp\ps{G}$-module which is finitely
generated over $\Zp\ps{H}$, then one has $\mu_G(M) =0$.

\item[$(c)$] Suppose that we are given a short exact sequence of finitely generated
$\Zp\ps{G}$ modules
\[ 0\lra M'\lra M\lra M'' \lra 0.\]
\begin{enumerate}
\item[$(i)$] One has $\mu_G(M) \leq \mu_G(M') + \mu_G(M'')$.
Moreover, if $M$, and hence also $M'$ and $M''$, is
$\Zp\ps{G}$-torsion, the inequality is an equality.

\item[$(ii)$] If $\mu_G(M'')=0$, then one has $\mu_G(M') =
\mu_G(M)$.

\item[$(iii)$] If $M'$ is finitely generated over $\Zp\ps{H}$, then one has $\mu_G(M)
=\mu_G(M'')$.

\end{enumerate}
\item[$(d)$] Suppose that we are given an exact sequence of finitely
generated $\Zp\ps{G}$-modules
\[ A\lra B \lra C\lra D\]
such that $A$ is finitely generated over $\Zp\ps{H}$ and
$\mu_G(D)=0$. Then one has the equality $\mu_G(B) = \mu_G(C)$.
\end{enumerate}
 \el

\bpf (a), (b) and (c)(i) are proven in \cite[Corollary 1.7]{Ho},
\cite[Lemma 2.7]{Ho} and \cite[Proposition 1.8]{Ho} respectively.
The remaining statements can be deduced from the previous statements
without too much difficulties. \epf

We make a note here mentioning that the conclusion of (c)(iii) is
false in general if one replaces the assumption ``$M'$ is finitely
generated over $\Zp\ps{H}$" by  ``$\mu_G(M')=0$". An example will be
to take $M' = M =\Zp\ps{G}$ and $M'' = \Zp\ps{G}/p$, and consider
the canonical exact sequence
\[ 0\lra \Zp\ps{G} \stackrel{p}{\lra} \Zp\ps{G} \lra \Zp\ps{G}/p\lra 0.\]
Clearly, we have $\mu_G(M')=0$ but $\mu_G(M) \neq \mu_G(M'')$.

\medskip
We end the section with the following useful relative formula for
the $\mu$-invariant.

\bl \label{mu relative}
 Let $G$ be a compact pro-$p$ $p$-adic Lie group
without $p$-torsion and let $N$ be a closed normal subgroup of $G$
such that $G/N$ has no $p$-torsion. Then for every finitely
generated torsion $\Zp\ps{G}$-module $M$, we have
\[  \mu_G(M) = \sum_{i\geq 0}(-1)^i \mu_{G/N}\big(H_i(N, M(p))\big).   \]\el

\bpf
 \[ \ba{rl}
   \mu_G(M) &= \displaystyle\sum_{i\geq 0}
   (-1)^i\mathrm{ord}_{p}\big(H_i(G,M(p))\big)\\
  & =  \displaystyle\sum_{i,j\geq 0}
   (-1)^{i+j}\mathrm{ord}_{p}\Big(H_i\big(G/N,
   H^j(N,M(p))\big)\Big)\\
  & =  \displaystyle\sum_{j\geq 0}
   (-1)^{j}  \displaystyle\sum_{i\geq 0}(-1)^i\mathrm{ord}_{p}\Big(H_i\big(G/N,
   H_j(N,M(p))\big)\Big)\\
   &  =  \displaystyle\sum_{j\geq 0}(-1)^j \mu_{G/N}\big(H_j(N, M(p))\big), \\
   \ea \]
where the first and fourth equality follow from Lemma \ref{mu
lemma}(a) (and noting that $H_j(N, M(p))$ is $p$-torsion), the third
is obvious and the second is a consequence of the following bounded
spectral sequence
 \[ H_i(G/N, H_j(N, M(p)))\Longrightarrow H_{i+j}(G,M(p)).\]
 \epf

\section{Selmer groups}

In this section, we recall the definition of the ($p$-primary)
Selmer group of an elliptic curve. As before, $p$ denote an odd
prime. Let $F$ be a number field, and let $E$ be an elliptic curve
defined over $F$. Let $v$ be a prime of $F$. For every finite
extension $L$ of $F$, we define
 \[ J_v(E/L) = \bigoplus_{w|v}H^1(L_w, E)(p),\]
where $w$ runs over the (finite) set of primes of $L$ above $v$. If
$\mathcal{L}$ is an infinite extension of $F$, we define
\[ J_v(E/\mathcal{L}) = \ilim_L J_v(E/L),\]
where the direct limit is taken over all finite extensions $L$ of
$F$ contained in $\mathcal{L}$. For any algebraic (possibly
infinite) extension $\mathcal{L}$ of $F$, the Selmer group of $E$
over $\mathcal{L}$ is defined to be
\[ S(E/\mathcal{L}) = \ker\Big(H^1(\mathcal{L}, E_{p^{\infty}})\lra \bigoplus_{v} J_v(E/\mathcal{L})
\Big), \] where $v$ runs through all the primes of $F$.

We say that $F_{\infty}$ is an \textit{admissible $p$-adic Lie
extension} of $F$ if (i) $\Gal(F_{\infty}/F)$ is a compact $p$-adic
Lie group, (ii) $F_{\infty}$ contains the cyclotomic $\Zp$-extension
$F^{\cyc}$ of $F$ and (iii) $F_{\infty}$ is unramified outside a
finite set of primes of $F$. Write $G = \Gal(F_{\infty}/F)$, $H =
\Gal(F_{\infty}/F^{\cyc})$ and $\Ga =\Gal(F^{\cyc}/F)$. Let $S$ be a
finite set of primes of $F$ which contains the primes above $p$, the
infinite primes, the primes at which $E$ has bad reduction and the
primes that are ramified in $F_{\infty}/F$. Denote $F_S$ to be the
maximal algebraic extension of $F$ unramified outside $S$. For each
algebraic (possibly infinite) extension $\mathcal{L}$ of $F$
contained in $F_S$, we write $G_S(\mathcal{L}) =
\Gal(F_S/\mathcal{L})$. The following alternative equivalent
description of the Selmer group of $E$
\[ S(E/\mathcal{L}) = \ker\Big(H^1(G_S(\mathcal{L}), E_{p^{\infty}})\stackrel{\la_S(F_{\infty})}{\lra} \bigoplus_{v\in S} J_v(E/\mathcal{L})
\Big)\] is well-known. We will denote $X(E/\mathcal{L})$ to be the
Pontryagin dual of $S(E/\mathcal{L})$.

\textbf{From now on, we will assume that for every prime $v$ of $F$
above $p$, our elliptic curve $E$ has either good ordinary reduction
or multiplicative reduction at $v$}. The following conjecture is
well-known.

\medskip \noindent \textbf{Conjecture.} \textit{$X(E/F^{\cyc})$ is a torsion
$\Zp\ps{\Ga}$-module.}

\medskip
The conjecture was first stated in \cite{Maz} for elliptic curves
that have good ordinary reduction at all primes of $F$ above $p$.
The general form we have here was stated in \cite{Sch} (see also
\cite{HO, OcV}). At present, the best result in support of the
conjecture is due to Kato \cite{K}, who has proven it when $F$ is
abelian over $\Q$ and $E$ is an elliptic curve defined over $\Q$
with good ordinary reduction at $p$.

One has a natural generalization of the above conjecture to
admissible $p$-adic extensions, namely that $X(E/F_{\infty})$ should
be $\Zp\ps{G}$-torsion for every admissible $p$-adic Lie extension
$F_{\infty}/F$. In this direction, this has been studied in
\cite{HV} when $G$ has dimension 2, and later in \cite{HO} when $G$
is a solvable uniform pro-$p$ group. However, such a generalization
of the conjecture does not allow one to define a suitable
characteristic element which is required in the formulation of the
main conjectures of non-commutative Iwasawa theory. To overcome this
difficulty, the following stronger conjecture, when $G$ has
dimension $> 1$, was introduced in \cite{CFKSV} for elliptic curve
having good ordinary reduction at all primes above $p$. In the case
when the elliptic curve has multiplicative reduction, this was
introduced in \cite{Lee} and was an important condition required in
order to apply Iwasawa-theoretical methods to the study of root
numbers (see also \cite{CFKS} for a similar study in the case when
the elliptic curve has good ordinary reduction at all primes above
$p$).

\medskip \noindent \textbf{$\M_H(G)$-Conjecture.} \textit{For every admissible $p$-adic Lie extension $F_{\infty}$
of $F$,  $X(E/F_{\infty})/X(E/F_{\infty})(p)$ is a finitely
generated $\Zp\ps{H}$-module.}

\medskip
The best evidence in support of the $\M_H(G)$-conjecture is in the
``$\mu=0$ situation'', i.e., the $\M_H(G)$-conjecture holds for
$X(E/F_{\infty})$ if $X(E/F^{\cyc})$ is finitely generated over
$\Zp$ and $F_{\infty}$ is a pro-$p$ extension of $F$ (see
\cite[Proposition 5.6]{CFKSV} or \cite[Theorem 2.1]{CS12}). In this
article, we will investigate the above conjecture from a general
point of view, that is namely the ``$\mu\neq 0$ situation". Our
first goal is to prove the following criterion for the conjecture to
hold. From now on, we write $X_f(E/F_{\infty}) =
X(E/F_{\infty})/X(E/F_{\infty})(p)$.

\bt \label{main} Let $p$ be an odd prime. Assume that $(i)$ $E$ has
either good ordinary reduction or multiplicative reduction at every
prime of $F$ above $p$, $(ii)$ $X(E/F^{\cyc})$ is
$\Zp\ps{\Ga}$-torsion, $(iii)$ $G$ is pro-$p$ without $p$-torsion,
$(iv)$ $H^2(G_S(F_{\infty}), E_{p^{\infty}})=0$ and $(v)$
$\la_S(F_{\infty})$ is surjective. Then $X_f(E/F_{\infty})$ is a
finitely generated $\Zp\ps{H}$-module if and only if
\[ \mu_G(X(E/F_{\infty})) = \mu_{\Ga}(X(E/F^{\cyc})) \]
and $H_i(H, X_f(E/F_{\infty}))$ is finitely generated over $\Zp$ for
every odd $i \leq \dim H-1$. \et

For the remainder of the section, we will compare our theorem with
existing known results, deferring the proof to Section \ref{Proof of
Theorem}. Before we begin the comparison, we first discuss the
relationship between the torsionness of the dual Selmer group and
the conditions: $H^2(G_S(F_{\infty}),E_{p^{\infty}}) =0$ and
surjectivity of $\la_S(F_{\infty})$. We begin with the following
proposition.

\bp \label{surjective implies torsion} Let $F_{\infty}/F$ be an
admissible $p$-adic Lie extension. If
$H^2(G_S(F_{\infty}),E_{p^{\infty}}) =0$ and $\la_S(F_{\infty})$ is
surjective, then $X(E/F_{\infty})$ is a torsion $\Zp\ps{G}$-module.
  \ep

\bpf This follows from a standard rank calculation noting the
formulas in \cite[Theorem 4.1]{OcV} and \cite[Proposition 7.4]{HV}.
 \epf

The next proposition is a partial converse to the preceding one.
Although we make use of the proposition mainly under case (i), we
feel it important to note down a few other cases. We mention that
the proof of the proposition under case (ii) is a standard
well-known argument which we have included for the reader's
convenience.

\bp \label{torsion implies surjective} Let $F_{\infty}/F$ be an
admissible $p$-adic Lie extension. Suppose that at least one of the
following statements holds.
\begin{enumerate}
 \item[$(i)$] $E_{p^{\infty}}$ is not rational over $F_{\infty}$.
 \item[$(ii)$]  For every $v\in S$, the decomposition group of $G$ at $v$ has
dimension $\geq 2$.
 \item[$(iii)$] $E$ has no additive reduction, $G$ is pro-$p$ and has no $p$-torsion, and the set $S$ consists
 precisely of the primes above
$p$, the infinite primes, the primes at which $E$ has bad reduction
and the primes that are ramified in $F_{\infty}/F$.
  \end{enumerate}
Then $X(E/F_{\infty})$ is a torsion $\Zp\ps{G}$-module if and only
if $H^2(G_S(F_{\infty}), E_{p^{\infty}}) =0$ and $\la_S(F_{\infty})$
is surjective. \ep

\bpf By Proposition \ref{surjective implies torsion}, it suffices to
show that if $X(E/F_{\infty})$ is a torsion $\Zp\ps{G}$-module, then
we have $H^2(G_S(F_{\infty}), E_{p^{\infty}}) =0$ and
$\la_S(F_{\infty})$ is surjective. We first consider the case that
$E_{p^{\infty}}$ is not rational over $F_{\infty}$. Then by
\cite[Proposition 10]{Z04}, we have that
$E(F_{\infty})_{p^{\infty}}$ is finite. The vanishing of
$H^2(G_S(F_{\infty}), E_{p^{\infty}})$ and the surjectivity of
$\la_S(F_{\infty})$ then follow from an application of \cite[Theorem
7.2]{HV}.

Now we will prove the above assertion under the assumption of (ii).
Since $X(E/F_{\infty})$ is a torsion $\Zp\ps{G}$-module by
hypothesis, we have that the dual fine Selmer group (see
\cite[Section 3]{CS05} for definition) is also a torsion
$\Zp\ps{G}$-module. By \cite[Lemma 3.1]{CS05}, this is equivalent to
$H^2(G_S(F_{\infty}), E_{p^{\infty}}) =0$. Now the
Cassels-Poitou-Tate sequence (for example, see \cite[1.7]{CS00})
gives an exact sequence
 \[ 0\lra S(E/F_{\infty})\lra H^1(G_S(F_{\infty}), E_{p^{\infty}})
 \stackrel{\la_S(F_{\infty})}{\lra} \bigoplus_{v\in S}J_v(F_{\infty}) \lra
 \big(\widehat{S}(E/F_{\infty})\big)^{\vee}
 \lra H^2(G_S(F_{\infty}), E_{p^{\infty}})\lra 0, \]
where $\widehat{S}(E/F_{\infty})$ is defined as the kernel of the
map
 \[ \plim_{L} H^1(G_S(L), T_p E) \lra \plim_{L}
 \bigoplus_{w|S}T_p H^1(L_w , E). \]
Since we have already shown that $H^2(G_S(F_{\infty}),
E_{p^{\infty}}) =0$, it follows from the Cassels-Poitou-Tate exact
sequence that $\coker \la_S(F_{\infty}) =
\widehat{S}(E/F_{\infty})^{\vee}$. One may now apply the formulas in
\cite[Theorem 4.1]{OcV} and \cite[Proposition 7.4]{HV} to conclude
that
 \[ \rank_{\Zp\ps{G}}H^1(G_S(F_{\infty}),E_{p^{\infty}})^{\vee} =
 [F:\Q] =\rank_{\Zp\ps{G}}\big(\bigoplus_{v\in S} J_v(F_{\infty})\big)^{\vee}. \]
Therefore, it follows that $X(E/F_{\infty})$ is a torsion
$\Zp\ps{G}$-module if and only if $\widehat{S}(E/F_{\infty})$ is a
torsion $\Zp\ps{G}$-module. On the other hand, the Poitou-Tate
sequence gives an exact sequence
 \[ H^2(G_S(F_{\infty}), E_{p^{\infty}})^{\vee} \lra   \plim_{L} H^1(G_S(L), T_p E)
  \stackrel{\varphi}{\lra} \plim_{L}
 \bigoplus_{w|S} H^1(L_w , T_p E). \] Since $H^2(G_S(F_{\infty}),
 E_{p^{\infty}})=0$, we have that $\varphi$ is injective. In particular, $\widehat{S}(E/F_{\infty})$ is contained
 in $\plim_{L}
 \bigoplus_{w|S} H^1(L_w , T_p E)$. Now by virtue of the assumption in (ii), we may apply \cite[Proposition 4.5]{OcV} to
 conclude that $\plim_{L} \bigoplus_{w|S} H^1(L_w , T_p E)$, and hence $\widehat{S}(E/F_{\infty})$, is $\Zp\ps{G}$-torsionfree.
 On other hand, we have shown above that
 $\widehat{S}(E/F_{\infty})$ is a torsion $\Zp\ps{G}$-module, and
 therefore, we must have $\widehat{S}(E/F_{\infty}) =0$. This gives
 the surjectivity of $\la_S(F_{\infty})$, as required.

Suppose that we are in the situation of (iii). Now if
$E_{p^{\infty}}$ is not rational over $F_{\infty}$, then we are
already done by (i). Therefore, we may assume that
$F(E_{p^{\infty}}) \sbs F_{\infty}$. It then follows from either
\cite[Lemma 2.8]{C} or \cite[Lemma 5.1]{CH} that the field
$F(E_{p^{\infty}})$, and hence $F_{\infty}$, has the property that
for each $v$ which is either above $p$ or is a bad reduction prime
for $E$ (noting that $E$ has no additive reduction), the dimension
of $G$ at $v$ is $\geq 2$. Now let $v$ be a prime of $F$ that is
ramified in $F_{\infty}$ and does not divide $p$. Since $v$ does not
divide $p$, it is unramified in $F^{\cyc}/F$. Therefore, every prime
$w$ of $F^{\cyc}$ above $v$ must ramify in $F_{\infty}/F^{\cyc}$.
Since $G$ has no $p$-torsion by one of the hypotheses in (iii), the
inertia group of $w$ in $F_{\infty}/F^{\cyc}$ must be infinite and
hence of dimension $\geq 1$. Adding this to the decomposition
component coming from $F^{\cyc}/F$ (since $v$ does not divide $p$),
it follows that the decomposition group of $v$ in $F_{\infty}/F$ has
dimension $\geq 2$. Thus, we are now in the situation of (ii), and
so we are done. \epf

We now note the following corollary and give some remarks on it.

\bc \label{torsion surjective corollary} Let $F_{\infty}/F$ be an
admissible $p$-adic Lie extension. If $X(E/L^{\cyc})$ is a torsion
$\Zp\ps{\Gal(L^{\cyc}/L)}$-module for every finite extension $L$ of
$F$ contained in $F_{\infty}$, then we have that
$H^2(G_S(F_{\infty}), E_{p^{\infty}}) =0$ and that
$\la_S(F_{\infty})$ is surjective. In particular, $X(E/F_{\infty})$
is a torsion $\Zp\ps{G}$-module. \ec

\bpf  Since $E_{p^{\infty}}$ is not rational over $L^{\cyc}$, it
follows from Proposition \ref{torsion implies surjective} that
$H^2(G_S(L^{\cyc}), E_{p^{\infty}}) = 0$ and $\la_S(L^{\cyc})$ is
surjective. Note that $H^2(G_S(F_{\infty}), E_{p^{\infty}}) =\ilim_L
H^2(G_S(L^{\cyc}), E_{p^{\infty}})$ and $\la_S(F_{\infty}) = \ilim_L
\la_S(L^{\cyc})$, where $L$ runs through all finite extensions of
$F$ contained in $F_{\infty}$. Therefore, we have that
$H^2(G_S(L^{\cyc}), E_{p^{\infty}}) = 0$ and $\la_S(L^{\cyc})$ is
surjective. The final assertion of the corollary then follows from
an application of Proposition \ref{surjective implies torsion}. \epf

It has been a long asked question (see \cite{CFKSV, CS12}) on
whether one can prove the $\M_H(G)$-conjecture under the hypothesis
that $X(E/L^{\cyc})$ is a torsion $\Zp\ps{\Gal(L^{\cyc}/L)}$-module
for every finite extension $L$ of $F$ contained in $F_{\infty}$.
Although we can do no better than Theorem \ref{main} and Corollary
\ref{torsion surjective corollary} in so far as showing the
$\M_H(G)$-conjecture, we are able to show that our dual Selmer group
exhibits ``$\M_H(G)$-like properties" (see Proposition \ref{Selmer
fg}). We finally remark that there is a related result of
Hachimori-Ochiai in the direction of Corollary \ref{torsion
surjective corollary}. Namely, their result \cite[Theorem 2.3]{HO}
states that if $X(E/F^{\cyc})$ is a torsion $\Zp\ps{\Ga}$-module and
$G$ is a uniform solvable pro-$p$ group, then $X(E/F_{\infty})$ is a
torsion $\Zp\ps{G}$-module.

\medskip We now describe how our Theorem \ref{main} compares with
existing results.

Let $H$ be pro-$p$ of dimension $1$ with no $p$-torsion. Such a
group $H$ is necessarily solvable. Then by \cite[Theorem 2.3]{HO}
(see also \cite[Theorem 2.8]{HV}), the $\Zp\ps{\Ga}$-torsionness of
$X(E/F^{\cyc})$ implies the $\Zp\ps{G}$-torsionness of
$X(E/F_{\infty})$. Therefore, if $E$ has either good ordinary
reduction or multiplicative reduction at every prime of $F$ above
$p$, $X(E/F^{\cyc})$ is $\Zp\ps{\Ga}$-torsion and $G$ is pro-$p$ of
dimension $\leq 2$ with no $p$-torsion, then we have that
$X_f(E/F_{\infty})$ is a finitely generated $\Zp\ps{H}$-module if
and only if
\[ \mu_G(X(E/F_{\infty})) = \mu_{\Ga}(X(E/F^{\cyc})). \]
This recovers \cite[Corollary 3.2]{CS12}.

We now consider the case that $H$ is of dimension either 2 or 3, and
has no $p$-torsion. Then under the same hypotheses as Theorem
\ref{main}, $X_f(E/F_{\infty})$ is a finitely generated
$\Zp\ps{H}$-module if and only if $H_1(H, X_f(E/F_{\infty}))$ is
finite and
\[ \mu_G(X(E/F_{\infty})) = \mu_{\Ga}(X(E/F^{\cyc})). \]
In this form, this criterion has been observed in \cite[Lemmas 5.3
and 5.4]{CFKSV} when $F_{\infty}$ is the field generated by all the
$p$-power division points of an elliptic curve without complex
multiplication.

\medskip
We conclude this section with two auxiliary results on the
$\mu$-invariant of the Selmer group. The first gives an inequality
under an extra condition of $F_{\infty}$.

\bp Let $p$ be an odd prime. Assume that $(i)$ $E$ has either good
ordinary reduction or multiplicative reduction at every prime of $F$
above $p$, $(ii)$ that $X(E/F^{\cyc})$ is $\Zp\ps{\Ga}$-torsion and
$(iii)$ that there is a finite family of closed normal subgroups
$H_i$ $(0\leq i\leq r)$ of $G$
 such that $1=H_0\sbs H_1 \sbs \cdots\sbs H_r =H$ and $H_i/H_{i-1}\cong
 \Zp$ for $1\leq i\leq r$. Then we have that
$H^2(G_S(F_{\infty}), E_{p^{\infty}})=0$ and $\la_S(F_{\infty})$ is
surjective. Furthermore, we have
\[ \mu_G(X(E/F_{\infty})) \leq \mu_{\Ga}(X(E/F^{\cyc})). \] \ep

\bpf This follows by applying Proposition \ref{main relative prop}
iteratively. \epf

The next result gives a necessary condition for $X_f(E/F_{\infty})$
to be finitely generated over $\Zp\ps{H}$.

\bp Retain the assumptions of Theorem \ref{main}. Suppose that
$X_f(E/F_{\infty})$ is a finitely generated $\Zp\ps{H}$-module. Then
for every finite extension $L$ of $F$ contained in $F_{\infty}$, we
have
 $$ \mu_{\Ga_L}\big(X(E/L^{\cyc})\big) = [L:F]\mu_{\Ga_F}\big(X(E/F^{\cyc})\big),$$
 where $\Ga_L = \Gal(L^{\cyc}/L)$.  \ep

\bpf Write $G_L= \Gal(F_{\infty}/L)$. Then we have
\[ \mu_{\Ga_L}\big(X(E/L^{\cyc})\big) = \mu_{G_L}(X(E/F_{\infty}))=
[L:F]\mu_{G}(X(E/F_{\infty})) = [L:F]\mu_{\Ga}(X(E/F^{\cyc})). \]
\epf

In the case when $G=\Zp^2$, the conclusion of the preceding
proposition turns out to be a sufficient condition for
$X_f(E/F_{\infty})$ to be finitely generated over $\Zp\ps{H}$ (cf.
\cite[Theorem 3.8]{CS12}). In view of this, a natural question will
be whether the conclusion of the preceding proposition is a
sufficient condition for $X_f(E/F_{\infty})$ to be a finitely
generated $\Zp\ps{H}$-module for a general $G$. We do not have an
answer to this at this point of writing.

\section{Proof of Theorem \ref{main}} \label{Proof of Theorem}

In this section, we will prove a result relating the quantities
$\mu_G(X(E/F_{\infty}))$ and $\mu_{\Ga}(X(E/F^{\cyc}))$. This
relationship will allow us to establish the required criterion in
Theorem \ref{main}. As before, $p$ will denote an odd prime number
and $E$ will denote an elliptic curve defined over $F$ with either
good ordinary reduction or multiplicative reduction at every prime
of $F$ above $p$.  We let $F_{\infty}$ denote a fixed admissible
$p$-adic Lie extension of $F$. Let $S$ be a finite set of primes of
$F$ which contains the primes above $p$, the infinite primes, the
primes at which $E$ has bad reduction and the primes that are
ramified in $F_{\infty}/F$. We continue to write $G =
\Gal(F_{\infty}/F)$, $H = \Gal(F_{\infty}/F^{\cyc})$ and $\Ga
=\Gal(F^{\cyc}/F)$. We split the section into three subsections,
where we will study the $H$-cohomology of $J_v(E/F_{\infty})$ and
$S(E/F_{\infty})$ in Subsection \ref{local cohomology} and
Subsection \ref{global cohomology} respectively. The calculations
done in these two subsections will be used in Subsection \ref{mu
relation} to give a relationship between the respective
$\mu$-invariants of $X(E/F_{\infty})$ and $X(E/F^{\cyc})$. We then
use this relationship to give a proof of Theorem \ref{main}. Again,
we like to remind the reader that our line of attack is deeply
inspired by the approach in \cite[Section 2]{CSS}.

\subsection{Local cohomology calculations} \label{local cohomology}

Before calculating the $H$-cohomology of the local terms, we give
another description of $J_v(\mathcal{L})$, where $\mathcal{L}$ is an
algebraic extension of $F^{\cyc}$. By \cite[P. 150]{CG}, for each
prime $v$ of $F$ above $p$, we have a short exact sequence
\[ 0 \lra C_v \lra E_{p^{\infty}} \lra D_v \lra 0\] of discrete
$\Gal(\bar{F}_v/F_v)$-modules which is characterized by the fact
that $C_v$ is divisible and that $D_v$ is the maximal quotient of
$E_{p^{\infty}}$ by a divisible subgroup such that the inertia group
acts on $D_v$ via a finite quotient. Since our elliptic curve $E$
has either good ordinary reduction or multiplicative reduction at
each prime of $F$ above $p$, we have that $C_v$ and $D_v$ are
divisible abelian groups of corank one. If fact, these groups can be
explicitly described when $E$ has either good ordinary reduction or
split multiplicative reduction at $v$. For instance, if $E$ has good
ordinary reduction at $v$, we may take $D_v$ to be
$\widetilde{E}_{v, p^{\infty}}$, where $\widetilde{E}_v$ is the
reduction of $E$ mod $v$, and take $C_v$ to be the kernel of the
natural surjection $E_{p^{\infty}}\tha \widetilde{E}_{v,
p^{\infty}}$. If $E$ has split multiplicative reduction at $v$, we
have $C_v = \mu_{p^{\infty}}$ and $D_v = \Qp/\Zp$ (cf.\ \cite[P.
69-70]{Gr99}).

Let $\mathcal{L}$ be an algebraic (possibly infinite) extension of
$F^{\cyc}$. For each non-archimedean prime $w$ of $\mathcal{L}$,
define $\mathcal{L}_w$ to be the union of the completions at $w$ of
the finite extensions of $F$ contained in $\mathcal{L}$. We can now
state the following lemma.

\bl \label{local} Let $\mathcal{L}$ be an algebraic extension of
$F^{\cyc}$ which is unramified outside a set of finite primes of
$F$. Then we have an isomorphism
\[ J_v(\mathcal{L}) \cong \begin{cases}   \ilim_\mathcal{L'}  \bigoplus_{w|v}H^1(\mathcal{L'}_w, D_v),& \mbox{if } v\mbox{ divides }p \\
      \ilim_\mathcal{L'}  \bigoplus_{w|v}H^1(\mathcal{L'}_w, E_{p^{\infty}}), & \mbox{if } v\mbox{ does not divides }p \end{cases}
 \] where the direct limit is taken over all finite extensions
$\mathcal{L}'$ of $F^{\cyc}$ contained in $\mathcal{L}$.\el

\bpf It suffices to prove the lemma in the case when $\mathcal{L}$
is a finite extension of $F^{\cyc}$. Now by definition, we have
\[J_v(\mathcal{L}) = \bigoplus_{w|v}H^1(\mathcal{L}_w, E)(p).\]
By the local Kummer theory of elliptic curves, we have an exact
sequence
\[ 0\lra E(\mathcal{L}_w)\ot\Qp/\Zp \stackrel{\kappa_{\mathcal{L}_w}}{\lra} H^1(\mathcal{L}_w, E_{p^{\infty}})
\lra H^1(\mathcal{L}_w, E)(p)\lra 0.\] Now if $v$ does not divide
$p$, then $E(\mathcal{L}_w)\ot\Qp/\Zp =0$. The required isomorphism
is then immediate in this case. Now suppose that $v$ divides $p$.
Since the profinite degree of $\mathcal{L}_w$ over $F_v$ is
divisible by $p^{\infty}$, it follows that
$\Gal(\bar{F_v}/\mathcal{L}_w)$ has $p$-cohomological dimension
$\leq 1$ (cf.\ \cite[Theorem 7.1.8(i)]{NSW}). Therefore, we obtain
an exact sequence
\[ H^1(\mathcal{L}_w, C_v) \stackrel{\varphi_{\mathcal{L}_w}}{\lra} H^1(\mathcal{L}_w, E_{p^{\infty}})
\lra H^1(\mathcal{L}_w, D_v)\lra 0.\] Now it follows from
\cite[Proposition 4.3]{CG} and \cite[P. 69-70]{Gr99} that $\im
\kappa_{\mathcal{L}_w} = \im \varphi_{\mathcal{L}_w}$. Hence we have
$ H^1(\mathcal{L}_w, E)(p) \cong H^1(\mathcal{L}_w, D_v)$, and this
gives the required conclusion for the case when $v$ divides $p$.
\epf

\bl \label{local 1}
 $H^i(H, \bigoplus_{v\in S}J_v(F_{\infty}))$ is cofinitely generated
 over $\Zp$ for every $i\geq 1$. Moreover, if $H$ has no $p$-torsion, then we have
 \[H^i(H, \bigoplus_{v\in S}J_v(F_{\infty})) = 0 \mbox{ for } i\geq d-1,\]
where $d$ is the dimension of $H$. \el

\bpf We denote $A_v$ to be $E_{p^{\infty}}$ or $D_v$ according as
$v$ does not or does divide $p$. By the Shapiro lemma, we have
 \[ H^i\big(H, \bigoplus_{v\in S}J_v(F_{\infty})\big) \cong \bigoplus_w
 H^i\big(H_w, H^1(F_{\infty, w}, A_v)\big), \]
where $w$ runs through the (finite) set of primes of $F^{\cyc}$
above $S$ and $H_w$ is the decomposition group in $H$ of some fixed
prime of $F_{\infty}$ lying above $w$. It then suffices to show that
$H^i\big(H_w, H^1(F_{\infty, w}, E_{p^{\infty}})(p)\big)$ is
cofinitely generated over $\Zp$ for every $i\geq 1$. For any
extension $\mathcal{L}$ of $F_w^{\cyc}$, the profinite degree of
$\mathcal{L}$ over $F_v$ is divisible by $p^{\infty}$, and
therefore, it follows that $\Gal(\bar{F_v}/\mathcal{L})$ has
$p$-cohomological dimension $\leq 1$ (cf.\ \cite[Theorem
7.1.8(i)]{NSW}). Hence the spectral sequence
\[ H^i\big(H_w, H^j(F_{\infty, w}, A_v)\big)\Longrightarrow H^{i+j}(F^{\cyc}_w, A_v)\]
degenerates to yield \[ H^i\big(H_w, H^1(F_{\infty, w},
A_v)\big)\cong H^{i+2}(H_w, A_v(F_{\infty, w})) \mbox{ for } i\geq
1.\] Since $H_w$ is a $p$-adic Lie group and $A_v(F_{\infty, w})$ is
cofinitely generated over $\Zp$, the latter cohomology group, and
hence $ H^i\big(H_w, H^1(F_{\infty, w}, A_v)\big)$, is cofinitely
generated over $\Zp$. This proves the first assertion. The second
assertion is also an immediate consequence from the above
isomorphism. \epf

We need another lemma on the cokernel of the restriction map of the
local cohomology groups.

\bl \label{local restrict}
 The cokernel of the restriction map $\bigoplus_{v\in S}J_v(F^{\cyc}) \stackrel{\ga}{\lra} \Big(\bigoplus_{v\in S}J_v(F_{\infty})\Big)^H$
 is cofinitely generated
 over $\Zp$.
\el

\bpf As before, we denote $A_v$ to be $E_{p^{\infty}}$ or $D_v$
according as $v$ does not or does divide $p$. By an application of
the Shapiro lemma and the Hochschild-Serre spectral sequence, we
have
 \[ \coker\ga \cong \bigoplus_{w}
 H^2\big(H_w, A_v(F_{\infty,w})\big), \]
where $w$ runs through the (finite) set of primes of $F^{\cyc}$
above $S$ and $H_w$ is the decomposition group in $H$ of some fixed
prime of $F_{\infty}$ lying above $w$. Clearly, each $H^2\big(H_w,
A_v(F_{\infty,w})\big)$ is cofinitely generated over $\Zp$ and since
the sum is a finite one, the lemma follows. \epf

\subsection{Global cohomology calculations} \label{global cohomology}

In this subsection, we will study the $H$-cohomology of the Selmer
group. We first record a lemma on the $H$-cohomology of
$H^1(G_S(F_{\infty}), E_{p^{\infty}})$.

\bl \label{1cohomology} If $H^2(G_S(F_{\infty}), E_{p^{\infty}}) =
H^2(G_S(F^{\cyc}),E_{p^{\infty}}) = 0$, then $H^i\big(H,
H^1(G_S(F_{\infty}), E_{p^{\infty}})\big)$ is cofinitely generated
over $\Zp$ for every $i\geq 1$. Moreover, if $H$ has no $p$-torsion,
then we have
 \[ H^i\big(H,
H^1(G_S(F_{\infty}), E_{p^{\infty}})\big) =0 \mbox{ for } i\geq
d-1,\] where $d$ is the dimension of $H$. \el

\bpf Under the given assumptions of the lemma, the spectral sequence
\[ H^i\big(H, H^j(G_S(F_{\infty}), E_{p^{\infty}})\big)\Longrightarrow
H^{i+j}(G_S(F^{\cyc}), E_{p^{\infty}})\] degenerates to yield \[
H^i\big(H, H^1(F_{\infty}, E_{p^{\infty}})\big)\cong H^{i+2}(H,
E_{p^{\infty}}(F_{\infty})) \mbox{ for } i\geq 1,\] and the latter
group is easily seen to be cofinitely generated over $\Zp$. The
second assertion is also immediate. \epf

\bp \label{cohomology of Selmer} Assume that hypotheses $(i)$,
$(ii)$, $(iv)$ and $(v)$ of Theorem \ref{main} hold. Then
$H^i\big(H, S(E/F_{\infty})\big)$ is cofinitely generated over $\Zp$
for every $i\geq 1$.  Moreover, if hypothesis (iii) of Theorem
\ref{main} holds, we then have
 \[ H^i(H, S(E/F_{\infty})) =0 \mbox{ for } i \geq \dim H.\]\ep

\bpf
 Consider the following commutative diagram
\[  \entrymodifiers={!! <0pt, .8ex>+} \SelectTips{eu}{}\xymatrix{
    0 \ar[r]^{} & S(E/F^{\cyc}) \ar[d] \ar[r] &  H^1(G_S(F^{\cyc}), E_{p^{\infty}}) \ar[d]
    \ar[r] & \bigoplus_{v\in S}J_v(F^{\cyc}) \ar[d]_{\ga} \ar[r] & 0 &\\
    0 \ar[r]^{} & S(E/F_{\infty})^H \ar[r]^{} & H^1(G_S(F_{\infty}), E_{p^{\infty}})^H \ar[r] & \
    \bigoplus_{v\in S}J_v(F_{\infty})^H \ar[r] &  H^1\big(H, S(E/F_{\infty})\big) \ar[r] & \cdots } \]
with exact rows, where the vertical maps are given by restriction.
Note that the top and bottom rows are exact by Proposition
\ref{torsion implies surjective}. To simplify notation, we write
$W_{\infty}= H^1(G_S(F_{\infty}), E_{p^{\infty}})$ and
 $J_{\infty} = \bigoplus_{v\in S}J_v(F_{\infty})$. By a diagram chasing argument,
 we have a long exact sequence
 \[ \coker\ga \lra H^1\big(H, S(E/F_{\infty})\big)
 \lra H^1(H, W_{\infty})
 \lra H^1(H, J_{\infty})\lra \cdots \]
 \[ \cdots\lra H^{i-1}(H, J_{\infty}) \lra H^i\big(H, S(E/F_{\infty})\big)
 \lra H^i(H, W_{\infty})
 \lra H^i(H, J_{\infty})\lra \cdots \]
It then follows from Lemmas \ref{local 1}, \ref{local restrict} and
\ref{1cohomology}, and the above long exact sequence that $H^i(H,
S(E/F_{\infty}))$ is cofinitely generated over $\Zp$ for every
$i\geq 1$. The second assertion follows from Lemmas \ref{local 1}
and \ref{1cohomology}. \epf

We record an immediate corollary.

\bc \label{cohomology Selmer corollary} Suppose that all the
hypotheses in Theorem \ref{main} hold. Denoting $d$ to be the
dimension of $H$, we have
 \[ H_i\big(H, X(E/F_{\infty})(p)\big) =0 \mbox{ for } i \geq d, \mbox{ and}\]
  \[ H_i\big(H, X_f(E/F_{\infty})\big) =0 \mbox{ for } i \geq d.\] \ec

\bpf The vanishing is clear for $i\geq d+1$. Now by Proposition
\ref{cohomology of Selmer}, we have $H_d\big(H,
X(E/F_{\infty})\big)=0$. Since $H_d(H, -)$ is left exact, and both
$X(E/F_{\infty})(p)$ and $X_f(E/F_{\infty})$ (= $p^r
X(E/F_{\infty})$ for some integer $r \geq  0$) are submodules of
$X(E/F_{\infty})$, we have the vanishing for $i=d$ too. \epf

\subsection{Relation between $\mu_G(X(E/F_{\infty}))$ and
$\mu_{\Ga}(X(E/F^{\cyc}))$} \label{mu relation}

We can now give the required relation (compare with
\cite[Proposition 2.13]{CSS}).

\bp \label{main prop} Let $p$ be an odd prime.  Assume that $(i)$
$E$ has either good ordinary reduction or multiplicative reduction
at every prime of $F$ above $p$, $(ii)$ $X(E/F^{\cyc})$ is
$\Zp\ps{\Ga}$-torsion, $(iii)$ $G$ has no $p$-torsion, $(iv)$
$H^2(G_S(F_{\infty}), E_{p^{\infty}})=0$ and $(v)$
$\la_S(F_{\infty})$ is surjective. Then we have that
$H_i(H,X_f(E/F_{\infty}))$ is a finitely generated
$\Zp\ps{\Ga}$-torsion module for all $i\geq 0$ and
\[ \mu_G(X(E/F_{\infty})) = \mu_{\Ga}(X(E/F^{\cyc})) +
\sum_{i=
0}^{d-1}(-1)^{i+1}\mu_{\Ga}\big(H_i(H,X_f(E/F_{\infty}))\big), \]
 where $d$ is the dimension of $H$.\ep

\bpf By Lemma \ref{mu relative}, we have
\[ \mu_G(X(E/F_{\infty})) = \displaystyle\sum_{i\geq
0}(-1)^i\mu_{\Ga}\big(H_i(H, X(E/F_{\infty})(p))\big). \]
 Write $X=
X(E/F_{\infty})$ and $X_f = X_f(E/F_{\infty})$. Taking $H$-homology
of the following short exact sequence
\[ 0\lra X(p)\lra X\lra X_f\lra 0,
\]
we obtain a long exact sequence
\[ \cdots \lra H_{i+1}(H, X)\lra H_{i+1}(H,
X_f) \lra H_i(H, X(p))\lra H_i(H, X)\lra \cdots\]
\[ \cdots\lra H_1(H, X) \lra H_{1}(H, X_f)\lra H_{0}(H,
X(p)) \lra H_0(H, X)\lra H_0(H, X_f)\lra 0.\]
 Clearly, $H_i(H, X(p))$ is a finitely generated
 $\Zp\ps{\Ga}$-torsion module. Combining this with Proposition \ref{cohomology of
 Selmer} and the above long exact sequence, we have that $H_i(H, X_f)$ is a finitely generated
 $\Zp\ps{\Ga}$-torsion module for $i\geq 1$. By a standard argument (for instance, see \cite[Lemma
2.4]{CS12}), the natural map
 $H_0(H, X(E/F_{\infty})) \lra X(E/F^{\cyc})$ has kernel and
 cokernel which are finitely generated over $\Zp$. Therefore, this
 implies that $H_0(H, X(E/F_{\infty}))$ is a $\Zp\ps{\Ga}$-torsion
 module and
 \[\mu_\Ga\big(H_0(H, X(E/F_{\infty}))\big) =
 \mu_{\Ga}\big(X(E/F^{\cyc})\big).\]
  Applying Lemma \ref{mu lemma}(b), (c) and  (d)
to the long exact sequence, we obtain
\[ \mu_{\Ga}\big(H_i(H, X(p))\big) = \mu_{\Ga}\big(H_{i+1}(H,
X_f)\big) \mbox{ for } i\geq 1, \mbox{ and } \]
\[ \mu_{\Ga}\big(H_0(H, X(p))\big) = \mu_{\Ga}\big(H_0(H,
X)\big) + \mu_{\Ga}\big(H_{0}(H, X_f)\big) -\mu_{\Ga}\big(H_{1}(H,
X_f)\big).\]  Combining these with Corollary \ref{cohomology Selmer
corollary}, we obtain the required
 equality.
 \epf

We can now prove Theorem \ref{main}.

\bpf[Proof of Theorem \ref{main}]
  To see that the ``if" direction holds, one observes that it follows from Lemma \ref{mu lemma}(b)
  and the equality in
  Proposition \ref{main prop} that $\mu_{\Ga}\big(H_0(H,X_f(E/F_{\infty}))\big)=0$. Since $H_0(H,X_f(E/F_{\infty}))$
  is a finitely generated $\Zp\ps{\Ga}$-module, this in turns implies that $H_0(H,X_f(E/F_{\infty}))$ is finitely generated over $\Zp$.
  Since $G$, and hence $H$, is pro-$p$, we may apply Nakayama Lemma to conclude that $X_f(E/F_{\infty})$ is finitely
  generated over $\Zp\ps{H}$. To see that the ``only if" direction
  holds, we see that if $X_f(E/F_{\infty})$ is finitely generated
  over $\Zp\ps{H}$, then $H_i(H, X_f(E/F_{\infty}))$ is finitely
  generated over $\Zp$ for all $i\geq 0$. By Lemma \ref{mu
  lemma}(b), these modules have trivial $\mu_{\Ga}$-invariant. Putting these into the equality
  in Proposition \ref{main prop}, we obtain $\mu_G(X(E/F_{\infty})) = \mu_{\Ga}(X(E/F^{\cyc}))$.
\epf

We mention that one can also prove an analogue of Proposition
\ref{main prop} replacing $F^{\cyc}$ by an intermediate admissible
subextension $F_{\infty}'$ of $F_{\infty}$. We will only state the
following special form of such a statement whose omitted proof is
similar to those in this section.

\bp \label{main relative prop} Let $p$ be an odd prime. Let
$F_{\infty}$ and $F_{\infty}'$ be two admissible $p$-adic extensions
of $F$ with $F_{\infty}'\subseteq F_{\infty}$ and
$N:=\Gal(F_{\infty}/F_{\infty}') \cong \Zp$. Assume that $(i)$ $E$
has either good ordinary reduction or multiplicative reduction at
every prime of $F$ above $p$, $(ii)$ $X(E/F^{\cyc})$ is
$\Zp\ps{\Ga}$-torsion, $(iii)$ $G$ and $G/N$ have no $p$-torsion,
$(iv)$ $H^2(G_S(F'_{\infty}), E_{p^{\infty}})=0$ and $(v)$
$\la_S(F'_{\infty})$ is surjective. Then we have that
$H^2(G_S(F_{\infty}), E_{p^{\infty}})=0$ and $\la_S(F_{\infty})$ is
surjective. Furthermore, we have
\[ \mu_G(X(E/F_{\infty})) = \mu_{G/N}(X(E/F_{\infty}'))-
\mu_{G/N}(X_f(E/F_{\infty})_N). \]
 \ep

\section{Akashi Series of Selmer Groups} \label{Akashi section}

Firstly, we record the following proposition.

\bp \label{Selmer fg} Assume that $(i)$ $E$ has either good ordinary
reduction or multiplicative reduction at every prime of $F$ above
$p$, $(ii)$ $X(E/F^{\cyc})$ is $\Zp\ps{\Ga}$-torsion, $(iii)$ $G$
has no $p$-torsion, $(iv)$ $H^2(G_S(F_{\infty}), E_{p^{\infty}})=0$
and $(v)$ $\la_S(F_{\infty})$ is surjective. Then the following
statements hold.

\begin{enumerate}
\item[$(a)$] $H_0(H, X(E/F_{\infty}))$ is $\Zp\ps{\Ga}$-torsion and
its $\mu_{\Ga}$-invariant is precisely the quantity
$\mu_{\Ga}\big(X(E/F^{\cyc})\big)$.
\item[$(b)$] $H_i(H, X(E/F_{\infty}))$ is finitely generated over
$\Zp$ for every $i\geq 1$.
\item[$(c)$] $H_i(H, X(E/F_{\infty}))=0$  for $i\geq \dim H$.
\end{enumerate}
\ep

\bpf (a) is shown in the proof of Proposition \ref{main prop}. (b)
and (c) are restatements of Proposition \ref{cohomology of Selmer}.
\epf

\br By the preceding proposition, we see that the Selmer groups
satisfy much stronger properties than a general module $M$ in
$\M_H(G)$ as in \cite[Lemma 3.1]{CFKSV}. Even so we are still not
able to show that $X(E/F_{\infty})$ lies in $\M_H(G)$ in general.\er

Now, assuming that all the hypotheses in Proposition \ref{Selmer fg}
are valid, we set $f_i$ to be the characteristic power series of the
module $H_i(H, X(E/F_{\infty}))$. Then the Akashi series of
$X(E/F_{\infty})$ is given by
\[ \prod_{i=0}^{d-1} f_i^{(-1)^i}, \]
where $d$ is the dimension of $H$. Note that it follows from
Proposition \ref{Selmer fg}(b) that $p$ does not divide $f_i$ for
every $i\geq 1$.

Combining results of Kato and Hachimori-Ochiai, we can find examples
where one can construct Akashi series of Selmer groups
unconditionally as follows: Assume for now that our elliptic curve
$E$ is defined over $\Q$ and has good ordinary reduction at $p$. For
a given finite abelian extension $F$ of $\Q$, Kato's result
\cite[Theorem 12.4]{K} tells us that $X(E/F^{\cyc})$ is a torsion
$\Zp\ps{\Ga}$-module. If $F_{\infty}$ is a solvable admissible
$p$-adic extension of $F$, it then follows from \cite[Theorem
2.3]{HO} that $X(E/F_{\infty})$ is a torsion $\Zp\ps{G}$-module. If
our elliptic curve satisfies one of the conditions in Proposition
\ref{torsion implies surjective} (for instance, $E$ has no complex
multiplication or no additive reduction), then hypotheses (iv) and
(v) of Proposition \ref{Selmer fg} are satisfied. Therefore, we can
define Akashi series for $X(E/F_{\infty})$ via the above discussion.

\medskip
We now mention how one can prove the results \cite[(1.1), Theorem
1.3 and Theorem 6.2]{Z11} by replacing the condition
``$X(E/F_{\infty})\in\M_H(G)$" with the (weaker) three conditions:
$X(E/F^{\cyc})$ is $\Zp\ps{\Ga}$-torsion, $H^2(G_S(F_{\infty}),
E_{p^{\infty}})=0$ and $\la_S(F_{\infty})$ is surjective. Going
through the arguments in \cite{Z11}, one sees that the condition
``$X(E/F_{\infty})\in\M_H(G)$" is used in the following ways.

(1) To define the Akashi series for $X(E/F_{\infty})$.

(2) To deduce that $X(E/F^{\cyc})$ is $\Zp\ps{\Ga}$-torsion which is
required in \cite[Proposition 5.1, Proposition 5.5]{Z11}.

(3) To yield $H^2(G_S(F^{\cyc}), E_{p^{\infty}})= 0$ which is
required in \cite[Lemma 5.2]{Z11}.

(4) To yield $H^2(G_S(F_{\infty}), E_{p^{\infty}})= 0$ which is
required in \cite[Lemma 5.2]{Z11}.

(5) To deduce that $\la_S(F_{\infty})$ is surjective which is
required in \cite[Formula (5.3), Lemma 5.7]{Z11}.

As seen above, the replaced conditions suffice for us to define the
Akashi series for $X(E/F_{\infty})$, and the necessary deductions in
(2)-(5) required for the argument in \cite{Z11} are consequences of
the replaced conditions by Propositions \ref{surjective implies
torsion} and \ref{torsion implies surjective}.

\medskip
We will now establish the following criterion result for the
vanishing of the Selmer group of an elliptic curve with good
ordinary reduction at all primes above $p$. Recall that a finitely
generated $\Zp\ps{G}$-module $M$ is said to be \textit{pseudo-null}
if $\Ext^i_{\Zp\ps{G}}(M, \Zp\ps{G}) =0$ for $i=0,1$.

\bt \label{Selmer vanishing} Assume that $(i)$ $E$ has good ordinary
reduction at every prime of $F$ above $p$, $(ii)$ $X(E/F^{\cyc})$ is
$\Zp\ps{\Ga}$-torsion, $(iii)$ $G$ is pro-$p$ with no $p$-torsion,
$(iv)$ $H^2(G_S(F_{\infty}), E_{p^{\infty}})=0$ and $(v)$
$\la_S(F_{\infty})$ is surjective. Furthermore, suppose that at
least one of the following statements holds.

\begin{enumerate}
\item[$(1)$] $X(E/F_{\infty})$ has no nonzero pseudo-null
$\Zp\ps{G}$-submodule.
\item[$(2)$] The number field $F$ is not totally real.
\end{enumerate}

Then $X(E/F_{\infty}) = 0$ if and only if its Akashi series is a
unit in $\Zp\ps{\Ga}$. In particular, if $X(E/F_{\infty}) \neq 0$,
then its Akashi series is not in $\Zp\ps{\Ga}^{\times}$. \et

We mention that condition (1) is satisfied in many cases (see
\cite[Theorem 3.2]{HO}, \cite[Theorem 2.6]{HV}, \cite[Theorem
6.5]{KT} and \cite[Theorem 5.1]{OcV02}). The above theorem under
condition (1) is probably known (for instance, see \cite[Proposition
A.9]{DD} or \cite[Proposition 4.12]{HV} for the case of a false Tate
extension), but nevertheless, we have included it for completeness
and also because the author could not find a reference for it in the
general case.

The proof of Theorem \ref{Selmer vanishing} will follow from a
series of lemmas. Before stating and proving these lemmas, we
introduce certain preliminary notation. Let $G$ be a compact pro-$p$
$p$-adic group without $p$-torsion, and let $H$ be a closed normal
subgroup of $G$ with $\Ga= G/H\cong \Zp$. For a given finitely
generated $\Zp\ps{G}$-module $M$, we denote $Ak(M)$ to be the Akashi
series of $M$ which is given by the alternating product of the
characteristic polynomials of its $H$-homology groups. Of course,
the Akashi series is only well-defined (up to a unit in
$\Zp\ps{\Ga}$) if all the $H$-homology groups of $M$ are
$\Zp\ps{\Ga}$-torsion. We can now state the following.

\bp \label{algebra prop} Let $G$ be a compact pro-$p$ $p$-adic group
without $p$-torsion, and let $H$ be a closed normal subgroup of $G$
with $G/H\cong \Zp$. Let $M$ be a finitely generated
$\Zp\ps{G}$-module which satisfies the following properties:
\begin{enumerate}
\item[$(i)$] $H_i(H,M)$ is a finitely generated $\Zp\ps{\Ga}$-torsion module for
every even $i$.
\item[$(ii)$] $H_i(H, M)$ is finitely generated over
$\Zp$ for every odd $i$.
\item[$(iii)$] $Ak(M)$ lies in $\Zp\ps{\Ga}^{\times}$.
\end{enumerate}
Then $M$ is a finitely generated torsion $\Zp\ps{H}$-module. In
particular, $M$ is a finitely generated pseudo-null
$\Zp\ps{G}$-module. \ep

We record the following corollary which will establish Theorem
\ref{Selmer vanishing} under condition (1).

\bc \label{algebra coro}
 Retain the notation and assumptions of the preceding
lemma. Suppose further that $M$ has no nonzero pseudo-null
$\Zp\ps{G}$-submodule. Then $Ak(M)$ is a unit in $\Zp\ps{\Ga}$ if
and only if $M =0$. \ec

\bpf
 One direction of the corollary is obvious. Conversely, suppose that $Ak(M)$ is
 a unit in $\Zp\ps{\Ga}$. Then Proposition \ref{algebra prop} tells
 us that $M$ is a finitely generated pseudo-null $\Zp\ps{G}$-module.
 Since $M$ has no nonzero pseudo-null
 $\Zp\ps{G}$-submodule, we must have $M=0$.
\epf

\br  Corollary \ref{algebra coro} is a refinement of an observation
made in \cite[P. 182]{CFKSV}. In the case when $G$ is abelian, the
converse of Proposition \ref{algebra prop} does hold (cf.
\cite[Lemma 4.4]{CSS}). When $G$ is not abelian, the converse  is
false in general (see \cite[Examples 3 and 4]{CSS}), although it is
known to hold in certain special cases (see \cite[Lemma 4.5]{CSS}
and \cite[Proposition 2.3]{Z11}). For another variant of Proposition
\ref{algebra prop}, we refer reader to Proposition \ref{algebra
complement} below. \er

From now on, we will identify $\Zp\ps{\Ga}\cong \Zp\ps{T}$ under a
choice of a generator of $\Ga$. A polynomial $T^n + c_{n-1}T^{n-1} +
\cdots +c_0$ in $\Zp[T]$ is said to be a Weierstrass polynomial if
$p$ divides $c_i$ for every $0\leq i \leq n-1$. We record the
following lemma.

\bl \label{algebra lemma}
 Let $f$ and $g$ be two Weierstrass polynomials, and let $a$ and $b$ be
 two non-negative integers. Then $p^af$ and $p^bg$ generate the same
 ideal in $\Zp\ps{T}$ if and only if $a=b$ and $(f) = (g)$.
 Furthermore, we have $\deg f =\deg g$.
\el

\bpf
 Since $p^af$ and $p^bg$ generate the same
 ideal in $\Zp\ps{T}$, we have a pseudo-isomorphism
 \[ \Zp\ps{T}/(p^a)\oplus \Zp\ps{T}/(f) \sim \Zp\ps{T}/(p^b)\oplus \Zp\ps{T}/(g). \]
This in turns implies that $a=b$ and $(f)=(g)$. Then we have
 \[ \deg f = \rank_{\Zp} \Zp\ps{T}/(f)= \rank_{\Zp} \Zp\ps{T}/(g)= \deg g,\]
 thus proving the lemma.  \epf

We can now give the proof of Lemma \ref{algebra prop}.

\bpf[Proof of Lemma \ref{algebra prop}]
 For each $i$, choose a Weierstrass polynomial $f_i$ such that
 $p^{a_i}f_i$ is a characteristic element of $H_i(H,M)$,
 where $a_i$ is the $\mu_{\Ga}$-invariant of $H_i(H,M)$. Write
 $$ a=\sum_{i~even} a_i, \quad f = \prod_{i~even} f_i \quad \mbox{and} \quad g=\prod_{i~odd}
 f_i.$$
Then we have $Ak(M) = p^a f/g$. Now suppose that $Ak(M)$ is a unit.
It then follows from Lemma \ref{algebra lemma} that $a=0$ and $\deg
f = \deg g$. In particular, we have $a_0 = 0$ which in turn implies
that $H_0(H,M)$ is a finitely generated $\Zp$-module. Since $G$ (and
hence $H$) is pro-$p$, we may apply Nakayama Lemma to conclude that
$M$ is finitely generated over $\Zp\ps{H}$. Therefore, it follows
that $H_i(H,M)$ is finitely generated over $\Zp$ and that $\deg f_i
= \rank_{\Zp} H_i(H,M)$. Then we have
 $$ \rank_{\Zp\ps{H}} M = \sum_{i\geq 0} (-1)^i\rank_{\Zp} H_i(H,M) =
 \sum_{i\geq 0} (-1)^i\deg f_i =
\deg f -\deg g =0,$$
 where the first equality follows from Howson's formula (cf. \cite[Theorem 1.1]{Ho}).
Therefore, we have established the first assertion of Proposition
\ref{algebra prop}. The second assertion follows from the first by a
well-known result of Venjakob (see \cite[Example 2.3 and Proposition
5.4]{V03}).\epf

It remains to show the validity of Theorem \ref{Selmer vanishing}
under condition (2), and this will follow from combining Lemma
\ref{algebra prop} with the following lemma. We note that this lemma
does not assume any additional condition on the structure of
$X(E/F_{\infty})$ (i.e., torsioness or $\M_H(G)$) other than being
nonzero.

\bl \label{algebra nonzero}  Assume that $E$ has good ordinary
reduction at every prime of $F$ above $p$, and assume that $F$ is
not totally real. If $X(E/F_{\infty})\neq 0$, then $X(E/F_{\infty})$
is not a finitely generated torsion $\Zp\ps{H}$-module.  \el

\bpf
 If $X(E/F_{\infty})$ is not finitely generated over $\Zp\ps{H}$, then we
are done. Therefore, we may assume that $X(E/F_{\infty})$ is
finitely generated over $\Zp\ps{H}$. In particular, this implies
that $X(E/L^{\cyc})$ is finitely generated over $\Zp$ for every
finite extension $L$ of $F$ contained in $F_{\infty}$. Note that
$X(E/F_{\infty}) = \plim_LX(E/L^{\cyc})$, where $L$ runs through all
finite extensions of $F$ contained in $F_{\infty}$. Since
$X(E/F_{\infty})\neq 0$, it follows that $X(E/L^{\cyc})\neq 0$ for
some $L$. Since $F$ is not totally real, so is $L$. Therefore, we
may apply \cite[Proposition 7.5]{Mat} to conclude that
$X(E/L^{\cyc})$ has positive $\Zp$-rank. By an application of
Hachimori's formula (see \cite[Theorem 5.4]{HS} or Proposition
\ref{Selmer general}), this in turn implies that $X(E/F_{\infty})$
has positive $\Zp\ps{H}$-rank. Hence we have established the lemma.
 \epf

The conclusion of Theorem \ref{Selmer vanishing} follows immediately
from Proposition \ref{algebra prop}, Corollary \ref{algebra coro}
and Lemma \ref{algebra nonzero}. We now proceed to establish an
analogous statement of Theorem \ref{Selmer vanishing} for the Artin
twists of the Selmer group under the assumption that the Selmer
group satisfies the $\M_H(G)$-conjecture. By abuse of notation, we
will denote $\M_H(G)$ to be the category of all finitely generated
$\Zp\ps{G}$-modules $M$ such that $M/M(p)$ is finitely generated
over $\Zp\ps{H}$. As a start, we record the following analog of
Proposition \ref{algebra prop}.

\bp \label{algebra complement} Let $G$ be a compact pro-$p$ $p$-adic
group without $p$-torsion, and let $H$ be a closed normal subgroup
of $G$ with $G/H\cong \Zp$. Let $M$ be a finitely generated
$\Zp\ps{G}$-module which lies in $\M_H(G)$. If $Ak(M)\in
\Zp\ps{\Ga}^{\times}$, then $M$ is a finitely generated pseudo-null
$\Zp\ps{G}$-module.

Furthermore, if one assumes that $M$ has no nontrivial pseudo-null
$\Zp\ps{G}$-submodule, then we have $Ak(M)\in \Zp\ps{\Ga}^{\times}$
if and only if $M = 0$. \ep

\bpf It suffices to show the first assertion. Write $M_f = M/M(p)$.
Then $Ak(M) = Ak(M(p))Ak(M_f)$. Note that $Ak(M(p)) = p^{\mu_G(M)}$
by \cite[Corollary 1.7]{Ho}. Since $M_f$ is finitely generated over
$\Zp\ps{H}$ by assumption, it follows that $Ak(M_f) = g/h$ for some
Weierstrass polynomials $g$ and $h$. Hence we have $Ak(M) =
p^{\mu_G(M)}g/h$. Since $Ak(M)\in \Zp\ps{\Ga}^{\times}$, it follows
from Lemma \ref{algebra lemma} that $\mu_G(M(p)) = \mu_G(M)=0$ and
$Ak(M_f)=g/h\in \Zp\ps{\Ga}^{\times}$. By \cite[Remark 3.33]{V02},
the equality $\mu_G(M(p))=0$ in turn implies that $M(p)$ is a
pseudo-null $\Zp\ps{G}$-module. On the other hand, since $Ak(M_f)\in
\Zp\ps{\Ga}^{\times}$ and $M_f$ is a finitely generated
$\Zp\ps{H}$-module, we may apply Proposition \ref{algebra prop} to
conclude that $M_f$ is a pseudo-null $\Zp\ps{G}$-module. Hence $M$
is a pseudo-null $\Zp\ps{G}$-module, since it is an extension of the
pseudo-null $\Zp\ps{G}$-modules $M(p)$ and $M_f$. \epf

Suppose that we are given an Artin representation $\rho:G\lra
GL_{d_{\rho}}(\Op_{\rho})$, where $\Op = \Op_{\rho}$ is the ring of
integers of some finite extension of $\Qp$. Denote $W_{\rho}$ to be
a free $\Op$-module of rank $d_{\rho}$ realizing $\rho$. If $M$ is a
$\Zp\ps{G}$-module, we define $\mathrm{tw}_{\rho}(M)$ to be the
$\Op_{\rho}$-module $W_{\rho}\ot_{\Zp}M$ with $G$ acting diagonally.

For an admissible $p$-adic extension $F_{\infty}$ with
$G=\Gal(F_{\infty}/F)$, we denote the twisted Selmer group by
$S(\mathrm{tw}_{\rho}(E)/F_{\infty})$ which is obtained by taking
$A=E_{p^{\infty}}\ot_{\Zp}W_{\rho}$ and $A_v = C_v\ot_{\Zp}W_{\rho}$
in the definition of the Greenberg Selmer group in Section
\ref{Complement Section} (see also \cite[P. 47-48]{CFKS}). We denote
$X(\mathrm{tw}_{\rho}(E)/F_{\infty})$ to be the Pontryagin dual of
$S(\mathrm{tw}_{\rho}(E)/F_{\infty})$. By \cite[Lemma 3.4]{CFKS}, we
have
$$\mathrm{tw}_{\hat{\rho}}X(E/F_{\infty}) =
X(\mathrm{tw}_{\rho}(E)/F_{\infty}),$$ where $\hat{\rho}$ is the
contragredient of $\rho$.

\bl \label{algebra nonzero twist}  Suppose that $E$ has good
ordinary reduction at every prime of $F$ above $p$ and suppose that
$F$ is not totally real. If $X(E/F_{\infty})\neq 0$, then
$X(\mathrm{tw}_{\rho}(E)/F_{\infty})$ is not a finitely generated
torsion $\Op\ps{H}$-module.  \el

\bpf
  We will prove the lemma by contradiction. Suppose that
$X(\mathrm{tw}_{\rho}(E)/F_{\infty})$ is a finitely generated
torsion $\Op\ps{H}$-module. Let $L$ be a finite extension of $F$
contained in $F_{\infty}$ such that $\Gal(F_{\infty}/L)$ is
contained in $\ker\rho$. Write $H_L = \Gal(F_{\infty}/L^{\cyc})$.
Since $H_L$ is a subgroup of $H$ with finite index and $\Op$ is a
finite free $\Zp$-algebra, we have that $\Op\ps{H}$ is a finite free
$\Zp\ps{H_L}$-algebra and

\[\ba{l} \Op\ps{H}\ot_{\Zp\ps{H_L}}\Hom_{\Zp\ps{H_L}}(X(\mathrm{tw}_{\rho}(E)/F_{\infty}),
\Zp\ps{H_L}) \vspace{0.1in} \\ \hspace{1in} \cong
\Hom_{\Op\ps{H}}(\Op\ps{H}\ot_{\Zp\ps{H_L}}X(\mathrm{tw}_{\rho}(E)/F_{\infty}),
\Op\ps{H}) \vspace{0.1in}\\ \hspace{1in} \cong
\Hom_{\Op\ps{H}}(X(\mathrm{tw}_{\rho}(E)/F_{\infty}),
\Op\ps{H})^{[H:H_L][\Op:\Zp]} \\ \hspace{1in}=0.
 \ea \]
This in turn implies that
$\Hom_{\Zp\ps{H_L}}(X(\mathrm{tw}_{\rho}(E)/F_{\infty}),
\Zp\ps{H_L})=0$. On the other hand, we have
\[\ba{l}0  =\Hom_{\Zp\ps{H_L}}(X(\mathrm{tw}_{\rho}(E)/F_{\infty}),
\Zp\ps{H_L}) \vspace{0.1in}\\ \hspace{0.1in}
 =\Hom_{\Zp\ps{H_L}}(X(E/F_{\infty})\ot_{\Zp}W_{\hat{\rho}}, \Zp\ps{H_L})
 \vspace{0.1in}\\ \hspace{0.1in}
\cong \Hom_{\Zp}\big(W_{\hat{\rho}},
\Hom_{\Zp\ps{H_L}}(X(E/F_{\infty}), \Zp\ps{H_L})\big).
  \ea \]
(The last isomorphism is the adjointness isomorphism which makes
sense here, since $H_L$ acts trivially on $W_{\hat{\rho}}$ by our
choice of $L$.) Since $W_{\hat{\rho}}$ is a free $\Zp$-module, it
follows from the above that
$$\Hom_{\Zp\ps{H_L}}(X(E/F_{\infty}),
\Zp\ps{H_L})=0.$$ This in turns implies that $X(E/F_{\infty})$ is a
torsion $\Zp\ps{H_L}$-module, contradicting Lemma \ref{algebra
nonzero}. Thus, we have proven our lemma. \epf

We can now prove the analogue of Theorem \ref{Selmer vanishing} for
the Artin twists of the Selmer group. This result may also be viewed
as a generalization of \cite[Corollary A.14]{DD}.

\bt \label{twist Selmer vanishing} Assume that $(i)$ $E$ has good
ordinary reduction at every prime of $F$ above $p$, $(ii)$ the
number field $F$ is not totally real, $(iii)$ $G$ is pro-$p$ and has
no $p$-torsion and $(iv)$ $X_f(E/F_{\infty})$ is finitely generated
over $\Zp\ps{H}$. Then the following statements are equivalent.

\begin{enumerate}

\item[$(a)$] $X(E/F_{\infty}) = 0$.

\item[$(b)$] $Ak(X(\mathrm{tw}_{\rho}(E)/F_{\infty}))$ is a unit in $\Op_{\rho}\ps{\Ga}$
for some Artin representation $\rho$ of $G$.

\item[$(c)$] $Ak(X(\mathrm{tw}_{\rho}(E)/F_{\infty}))$ is a unit in $\Op_{\rho}\ps{\Ga}$
for every Artin representation $\rho$ of $G$. \end{enumerate} \et

\bpf Clearly, one has the implications
(a)$\Rightarrow$(c)$\Rightarrow$(b). Now suppose that (b) holds. For
convenience, we will write $\Op=\Op_{\rho}$. Since
$X_f(E/F_{\infty})$ is finitely generated over $\Zp\ps{H}$, it
follows from \cite[Lemma 3.2]{CFKSV} that
$X_f(\mathrm{tw}_{\rho}(E)/F_{\infty})$ is a finitely generated
$\Op\ps{H}$-module. Therefore, we may apply an $\Op$-analogue of
Theorem \ref{main} to conclude that
\[\mu_G \big(
X(\mathrm{tw}_{\rho}(E)/F_{\infty})\big) = \mu_{\Ga} \big(
X(\mathrm{tw}_{\rho}(E)/F^{\cyc})\big). \] We then apply an
$\Op$-analogue of Proposition \ref{algebra complement} to conclude
that that $X(\mathrm{tw}_{\rho}(E)/F_{\infty})$ is a finitely
generated pseudo-null $\Op\ps{G}$-module. In particularly, it
follows from the proof of Proposition \ref{algebra complement} that
one has $\mu_G \big(X(\mathrm{tw}_{\rho}(E)/F_{\infty})\big) = 0$.
Combining this with the above, we obtain $\mu_{\Ga} \big(
X(\mathrm{tw}_{\rho}(E)/F^{\cyc})\big)=0$. By the structure theory
of finitely generated $\Op\ps{\Ga}$-modules, this in turn implies
that $X(\mathrm{tw}_{\rho}(E)/F^{\cyc})$ is finitely generated over
$\Op$. It then follows from an application of Nakayama Lemma that
$X(\mathrm{tw}_{\rho}(E)/F_{\infty})$ is finitely generated over
$\Op\ps{H}$. Since $X(\mathrm{tw}_{\rho}(E)/F_{\infty})$ is also a
pseudo-null $\Op\ps{G}$-module, we can apply a well-known result of
Venjakob (see \cite[Example 2.3 and Proposition 5.4]{V03}) to
conclude that $X(\mathrm{tw}_{\rho}(E)/F_{\infty})$ is a finitely
generated torsion $\Op\ps{H}$-module. By Lemma \ref{algebra nonzero
twist}, this in turn implies that $X(E/F_{\infty}) = 0$. \epf

\medskip
We end the section considering an analogue of Proposition
\ref{Selmer fg} for the Galois group of the maximal abelian
extension of $F_{\infty}$ unramified outside $S$. The statement may
have been known among the experts but does not seem to be written
down anywhere.

\bp \label{Galois fg} Assume that $F_{\infty}$ is an admissible
extension of $F$ such that $G = \Gal(F_{\infty}/F)$ has no
$p$-torsion. Let $S$ be a finite set of primes of $F$ containing
those above $p$, the infinite primes and the primes that ramify in
$F_{\infty}/F$. Write $X_S(F_{\infty}) = G_S(F_{\infty})^{ab}(p)$.
Then the following statements hold.

\begin{enumerate}
\item[$(a)$] $H_0(H, X_S(F_{\infty}))$ has $\Zp\ps{\Ga}$-rank $r_2(F)$ and
its $\mu_{\Ga}$-invariant is precisely the quantity
$\mu_{\Ga}\big(X_S(F^{\cyc})\big)$.
\item[$(b)$] $H_i(H, X_S(F_{\infty}))$ is finitely generated over
$\Zp$ for $i\geq 1$.
\item[$(c)$] $H_i(H, X_S(F_{\infty}))=0$  for $i\geq \dim G-1$.
\end{enumerate}
\ep

\bpf
 It is a fundamental knowledge that $H^2(G_{S}(F^{\cyc}),
 \Qp/\Zp)= 0$ and $H^2(G_{S}(F_{\infty}),
 \Qp/\Zp)= 0$. Therefore, the spectral sequence
 \[ H^i\big(H, H^j(G_S(F_{\infty}), \Qp/\Zp)\big) \Longrightarrow
 H^{i+j}(G_S(F^{\cyc}), \Qp/\Zp) \] degenerates to give an exact
 sequence
 \[ 0 \lra H^1(H,\Qp/\Zp)\lra H^1(G_S(F^{\cyc}), \Qp/\Zp) \lra H^1(G_S(F_{\infty}),
 \Qp/\Zp)^H \lra H^2(H, \Qp/\Zp) \] and
 an isomorphism
 \[ H^i\big(H, H^1(G_S(F_{\infty}), \Qp/\Zp)\big) \cong H^{i+2}(H,
 \Qp/\Zp) \mbox{ for each } i\geq 1.\]
All the statements in the proposition now follow immediately from
the facts that $H^i(H, \Qp/\Zp)$ is cofinitely generated over $\Zp$
and that $H^1(G_S(F_{\infty}), \Qp/\Zp)^{\vee} \cong
X_S(F_{\infty})$. \epf

\section{On the integrality of the
Akashi series} \label{Integrality of Akashi series}

In this section, we discuss certain integrality properties of the
Akashi series of Selmer groups. We will also establish
\cite[Conjecture 4.8 (Case 4)]{CFKSV} for the characteristic
elements attached to certain classes of modules in $\M_H(G)$.

As before, we say that $F_{\infty}$ is an admissible $p$-adic Lie
extension of $F$ if (i) $\Gal(F_{\infty}/F)$ is a compact $p$-adic
Lie group, (ii) $F_{\infty}$ contains the cyclotomic $\Zp$-extension
$F^{\cyc}$ of $F$ and (iii) $F_{\infty}$ is unramified outside a
finite set of primes of $F$. We continue to write $G =
\Gal(F_{\infty}/F)$, $H = \Gal(F_{\infty}/F^{\cyc})$ and $\Ga
=\Gal(F^{\cyc}/F)$.

Let $\rho:G\lra GL_{d_{\rho}}(\Op_{\rho})$ be an Artin
representation, where $\Op= \Op_{\rho}$ is the ring of integers of
some finite extension of $\Qp$. As in the previous section, the
twisted Selmer group $S(\mathrm{tw}_{\rho}(E)/F_{\infty})$ is
defined by taking $A=E_{p^{\infty}}\ot_{\Zp}W_{\rho}$ and $A_v =
C_v\ot_{\Zp}W_{\rho}$ in the definition of the Greenberg Selmer
group in Section \ref{Complement Section} (see also \cite[P.
47-48]{CFKS}). We denote by $\la_{\rho}(F_{\infty})$ the
localization map and by $X(\mathrm{tw}_{\rho}(E)/F_{\infty})$ the
Pontryagin dual of $S(\mathrm{tw}_{\rho}(E)/F_{\infty})$. We can now
state the following proposition which slightly refines \cite[Lemma
3.9]{CFKSV} for $X(E/F_{\infty})$.

\bp \label{integral Akashi} Assume that $(i)$ $E$ has good ordinary
reduction at every prime of $F$ above $p$, $(ii)$
$G=\Gal(F_{\infty}/F)$ is pro-$p$ and has no $p$-torsion, $(iii)$
$X_f(E/F_{\infty})$ is finitely generated over $\Zp\ps{H}$ and
$(iv)$ $H_i(H', X(E/F_{\infty}))$ is finite for every open normal
subgroup $H'$ of $H$ and $i\geq 1$.

Then $Ak(\mathrm{tw}_{\hat{\rho}}(X(E/F_{\infty}))) \in
\Op\ps{\Ga}/\Op\ps{\Ga}^{\times}$ for every Artin representation
$\rho$ of $G$. \ep

\bpf
 Since $X_f(E/F_{\infty})$ is finitely generated over $\Zp\ps{H}$, it
follows from \cite[Lemma 3.2]{CFKSV} and \cite[Lemma 3.4]{CFKS} that
$X_f(\mathrm{tw}_{\rho}(E)/F_{\infty})$ is a finitely generated
$\Op\ps{H}$-module. By an application of \cite[Proposition
2.5]{CS12}, this in turn implies that
$X(\mathrm{tw}_{\rho}(E)/L^{\cyc})$ is a finitely generated torsion
$\Op\ps{\Gal(L^{\cyc}/L)}$-module for every finite extension $L$ of
$F$ contained in $F_{\infty}$.

We claim that $(E[p^{\infty}]\ot_{\Zp}W_{\rho})(L^{\cyc})$ is
finite. Suppose for now that this claim holds. Then by a similar
argument to that of Corollary \ref{torsion surjective corollary}, we
have that $H^2(G_S(L^{\cyc}),E[p^{\infty}]\ot_{\Zp}W_{\hat{\rho}})
=0$ and the localization map $\la_{\rho}(L^{\cyc})$ is surjective
for all $L$ (and in particular for $F$), and that
$H^2(G_S(F_{\infty}),E[p^{\infty}]\ot_{\Zp}W_{\rho}) =0$ and the
localization map $\la_{\rho}(F_{\infty})$ is surjective. Therefore,
by an entirely parallel argument to that of Proposition \ref{Selmer
fg}, we conclude that
$H_i(H,\mathrm{tw}_{\hat{\rho}}(X(E/F_{\infty})))$ is finitely
generated over $\Op$ for every $i\geq 1$. On the other hand, by
virtue of the assumptions (iii) and (iv), it follows from the proof
of \cite[Lemma 3.9]{CFKSV} that, for every $i\geq 1$,
$H_i(H,\mathrm{tw}_{\hat{\rho}}(X(E/F_{\infty})))$ is annihilated by
some power of $p$. Hence we conclude that
$H_i(H,\mathrm{tw}_{\hat{\rho}}(X(E/F_{\infty})))$ is finite for
every $i\geq 1$. Therefore, the Akashi series of
$\mathrm{tw}_{\hat{\rho}}(X(E/F_{\infty}))$ is precisely the
$\Zp\ps{\Ga}$-characteristic power series of $H_0(H,
\mathrm{tw}_{\hat{\rho}}(X(E/F_{\infty})))$, and so, necessarily
lies in $\Op\ps{\Ga}/\Op\ps{\Ga}^{\times}$.

It remains to verify our claim. Since $\rho$ is an Artin
representation, we may find some finite extension $L_0$ of $L$
contained in $F_{\infty}$ such that $\rho$ factors though
$\Gal(L_0/L)$. Then we have
$$(E_{p^{\infty}}\ot_{\Zp}W_{\rho})(L_0^{\cyc})= E_{p^{\infty}}(L_0^{\cyc})^{[\Op:\Zp]},$$
where the latter is finite by \cite[Theorem 4.3]{Wi}. Thus, it
follows that $(E_{p^{\infty}}\ot_{\Zp}W_{\rho})(L^{\cyc})$ is also
finite. \epf

To prove the next set of results, we need to introduce some further
notion and notation. Let
\[ \Si = \{\,s\in \Zp\ps{G}~ | ~\Zp\ps{G}/\Zp\ps{G}s
~\mbox{is a finitely generated $\Zp\ps{H}$-module} \}.\]
 By
\cite[Theorem 2.4]{CFKSV}, $\Si$ is a left and right Ore set
consisting of non-zero divisors in $\Zp\ps{G}$. Set $\Si^* =
\cup_{n\geq 0}p^n\Si$. It follows from \cite[Proposition 2.3]{CFKSV}
that a finitely generated $\Zp\ps{G}$-module $M$ is annihilated by
$\Si^*$ if and only if $M/M(p)$ is finitely generated over
$\Zp\ps{H}$. Therefore, the category $\M_H(G)$ can also be thought
as the category of all finitely generated $\Zp\ps{G}$-modules which
are $\Si^*$-torsion. By the discussion in \cite[Section 3]{CFKSV},
we have the following exact sequence
\[ K_1(\Zp\ps{G}) \lra K_1(\Zp\ps{G}_{\Si^*})\stackrel{\partial_G}{\lra} K_0(\M_H(G))\lra 0 \]
of $K$-groups. For each $M$ in $\M_H(G)$, we define a characteristic
element for $M$ to be any element $\xi_M$ in
$K_1(\Zp\ps{G}_{\Si^*})$ with the property that
\[\partial_G(\xi_M) = [M].\]

Let $\rho:G\lra GL_{d_{\rho}}(\Op_{\rho})$ denote a continuous group
representation (not necessarily an Artin representation), where
$\Op= \Op_{\rho}$ is the ring of integers of some finite extension
of $\Qp$. For $g\in G$, we write $\bar{g}$ for its image in
$\Ga=G/H$. We define a continuous group homomorphism
 \[G \lra M_d(\Op)\ot_{\Zp}\Zp\ps{\Ga}, \quad g\mapsto \rho(g)\ot
 \bar{g}. \] By \cite[Lemma 3.3]{CFKSV}, this in turn induces a map
 \[ \Phi_{\rho}': K_1(\Zp\ps{G}_{\Si^*})\lra Q_{\Op}(\Ga), \]
where $Q_{\Op}(\Ga)$ is the field of fraction of $\Op\ps{\Ga}$. Let
$\varphi: \Op\ps{\Ga}\lra \Op$ be the augmentation map and denote
its kernel by $\mathfrak{p}$. One can extend $\varphi$ to a map
$\varphi : \Op\ps{\Ga}_{\mathfrak{p}}\lra K$, where $K$ is the field
of fraction of $\Op$. Let $\xi$ be an arbitrary element in
$K_1(\Zp\ps{G}_{\Si^*})$. If $\Phi_{\rho}'(\xi)\in
\Op\ps{\Ga}_{\mathfrak{p}}$, we define $\xi(\rho)$ to be
$\varphi(\Phi_{\rho}'(\xi))$. If $\Phi_{\rho}'(\xi)\notin
\Op\ps{\Ga}_{\mathfrak{p}}$, we set $\xi(\rho)$ to be $\infty$.

We will write $\al$ for the natural map
\[\Zp\ps{G}_{\Si^*}^{\times}\lra K_1(\Zp\ps{G}_{\Si^*}).\] We can now state our
results which will prove \cite[Conjecture 4.8 Case 4]{CFKSV} for the
characteristic elements attached to certain modules.

\bp \label{conjecture 4.8} Let $M$ be an object in $\M_H(G)$. Assume
that $M$ has no nonzero pseudo-null $\Zp\ps{G}$-submodule. Let
$\xi_M$ be a characteristic element of $M$. Then the following
statements are equivalent.

\begin{enumerate}
\item[$(a)$] $\xi_M \in \al(\Zp\ps{G}^{\times})$, where $\al$ is
the map $\Zp\ps{G}_{\Si^*}^{\times}\lra K_1(\Zp\ps{G}_{\Si^*})$.

\item[$(b)$] $\xi_M(\rho)$ is finite and lies in $\Op_{\rho}$ for every continuous
group representation $\rho$ of $G$.

\item[$(c)$] $\Phi_{\rho}'(\xi_M) \in \Op_{\rho}\ps{\Ga}^{\times}$ for every continuous
group representation $\rho$ of $G$.

\item[$(d)$] $\Phi_{\rho}'(\xi_M) \in \Op_{\rho}\ps{\Ga}^{\times}$ for
 every Artin representation $\rho$ of $G$.
\end{enumerate}
\ep

\bpf
 By \cite[Lemma 4.9]{CFKSV}, we have the implications
(a)$\Rightarrow$(b)$\Leftrightarrow$(c)$\Rightarrow$(d). It remains
to show (d)$\Rightarrow$(a). Taking $\rho:G\lra \Zp^{\times}$ to be
the trivial representation, one has
\[ \Phi_{\rho}(\xi_M) =
Ak(M)~\mathrm{mod}\,\Zp\ps{\Ga}^{\times}\] for this particular
$\rho$ by \cite[Lemma 3.7]{CFKSV}. Applying (d), we have that
$Ak(M)\in\Zp\ps{\Ga}^{\times}$. By Corollary \ref{algebra coro},
this in turn implies that $\partial_G(\xi_M)=0$. It then follows
from the above exact sequence of $K$-groups that there exists an
element in $K_1(\Zp\ps{G})$ which maps to $\xi_M$. On the other
hand, it is well-known that $\Zp\ps{G}^{\times}$ maps onto
$K_1(\Zp\ps{G})$, thus yielding (a). \epf

In particular, if $X(E/F_{\infty})$ lies in $\M_H(G)$ and has no
nonzero pseudo-null $\Zp\ps{G}$-module, the above proposition
applies to the characteristic elements of $X(E/F_{\infty})$. In the
case when we do not know whether $X(E/F_{\infty})$ has no nonzero
pseudo-null $\Zp\ps{G}$-module, we have the following.

\bp \label{conjecture 4.8(a)} Assume that $(i)$ $E$ has good
ordinary reduction at every prime of $F$ above $p$, $(ii)$ the
number field $F$ is not totally real, $(iii)$ $G$ is pro-$p$ and has
no $p$-torsion and $(iv)$ $X(E/F_{\infty})\in \M_H(G)$. Let $\xi_E$
be a characteristic element of $X(E/F_{\infty})$. Then the following
statements are equivalent.

\begin{enumerate}
\item[$(a)$] $\xi_E \in \al(\Zp\ps{G}^{\times})$, where $\al$ is
the map $\Zp\ps{G}_{\Si^*}^{\times}\lra K_1(\Zp\ps{G}_{\Si^*})$.

\item[$(b)$] $\xi_E(\rho)$ is finite and lies in $\Op_{\rho}$ for every continuous
group representation $\rho$ of $G$.

\item[$(c)$] $\Phi_{\rho}'(\xi_E) \in \Op_{\rho}\ps{\Ga}^{\times}$ for every continuous
group representation $\rho$ of $G$.

\item[$(d)$] $\Phi_{\rho}'(\xi_E) \in \Op_{\rho}\ps{\Ga}^{\times}$ for every Artin
 representation $\rho$ of $G$.
\end{enumerate}
\ep

\bpf
 The proof is similar to that of the preceding proposition, where one makes use of
Theorem \ref{twist Selmer vanishing} in place of Corollary
\ref{algebra coro}. \epf

We end the section mentioning that under the assumptions of
Proposition \ref{integral Akashi}, one has $\Phi_{\rho}'(\xi_E) \in
\Op_{\rho}\ps{\Ga}$ which in turn implies that assertion (d) of
\cite[Conjecture 4.8 Case 2]{CFKSV} is satisfied. However, the
author is not able to prove assertions (a), (b) and (c) of
\cite[Conjecture 4.8 Case 2]{CFKSV} (or even case 1 of the said
conjecture) for $\xi_E$.

\section{Complement: Other Ordinary Representations}
\label{Complement Section}
 As seen in the previous section,
most of the argument in this article can be applied to other
situations. Indeed, this is the case, and we will briefly describe
here a general context that one may apply the argument to.

Denote $\Op$ to be the ring of integers of some finite extension $K$
of $\Qp$. Fix a local parameter $\pi$ for $\Op$. The material in
Section \ref{algebra} can be extended over for modules over
$\Op\ps{G}$. One simply replaces any occurrences of ``$p$" by
``$\pi$", ``$\Zp$" by ``$\Op$" and ``$\mathbb{F}_p$" by
``$\Op/\pi\Op$".

Let $A$ be a cofinitely generated $\Op$-module with a continuous,
$\Op$-linear $\Gal(\bar{F}/F)$-action which is unramified outside a
finite set of primes of $F$. Assume further that $A^{\vee}$ has
$\Op$-rank $2d$, and for each $v$ above $p$, there is a
$\Op$-submodule $A_v$ of $A$ invariant under the action of
$\Gal(\bar{F}_v/F_v)$ and such that $(A_v)^{\vee}$ has $\Op$-rank
$d$. Write $D_v = A/A_v$. Note that $(D_v)^{\vee}$ has $\Op$-rank
$d$.

For an algebraic (possibly infinite) extension $\mathcal{L}$ of
$F^{\cyc}$, we define the Greenberg Selmer group of $A$ over
$\mathcal{L}$ by
\[ S(A/\mathcal{L}) = \ker\Big(H^1(\mathcal{L}, A)\lra \prod_{w\nmid p} H^1(\mathcal{L}, A)\times \prod_{w\mid p}H^1(\mathcal{L},
D_w)\Big), \] where we set $D_w = D_v$ whenever $w$ divides $v$. Now
for every finite extension $\mathcal{L}$ of $F^{\cyc}$, we define
$J_v(A/\mathcal{L})$ to be
 \[  \bigoplus_{w|v}H^1(\mathcal{L}_w, A)~\mbox{or}~ \bigoplus_{w|v}H^1(\mathcal{L}_w, D_w)\]
according as $v$ does not or does divide $p$. In the case that
$\mathcal{L}$ is an infinite extension of $F^{\cyc}$, we will define
\[ J_v(A/\mathcal{L}) = \ilim_L J_v(A/\mathcal{L}'),\]
where the direct limit is taken over all finite extensions
$\mathcal{L}'$ of $F^{\cyc}$ contained in $\mathcal{L}$. Let $S$ be
a finite set of primes of $F$ which contains the primes above $p$,
the ramified primes of $A$ and the infinite primes. It then follows
from a standard argument (see \cite[Corollary 3.2]{CS12}) that for
every extension $\mathcal{L}$ of $F^{\cyc}$ contained in $F_S$, one
has the following exact sequence
\[ 0 \lra S(A/\mathcal{L}) \lra H^1(G_S(\mathcal{L}), A)\stackrel{\la_S(\mathcal{L})}{\lra} \bigoplus_{v\in S} J_v(A/\mathcal{L}).\]

We will write $X(A/\mathcal{L})= S(A/\mathcal{L})^{\vee}$. Then one
can prove results analogous to Theorem \ref{main} and Proposition
\ref{Selmer fg} for $X(A/\mathcal{L})$. Note that we may have to
replace the condition ``$X(A/F^{\cyc})$ is $\Op\ps{\Ga}$-torsion" by
the conditions ``$H^2(G_S(F^{\cyc}), A)=0$" and ``$\la_S(F^{\cyc})$
is surjective". This is because we do not have an analog statement
as in Proposition \ref{torsion implies surjective}(i) for a general
$A$. The point is that we may not have the finiteness of
$A(F^{\cyc})$. However, for many interesting representations which
we shall see below, this does hold.

Typical examples of the above (besides the case of an elliptic curve
and its Artin twist as considered in the main body of the paper) are
abelian varieties $A$ with good ordinary reduction at all primes
above $p$. We note that the finiteness of $A(F^{\cyc})$ follows from
\cite[Theorem 4.3]{Wi}. The nonexistence of nontrivial pseudo-null
submodules for the dual Selmer groups was established in certain
cases \cite{Oc}, although they contain less situation of the
$p$-adic extensions than the case of an elliptic curve. Despite
this, if one assumes that the base field $F$ is not totally real,
one can still make use of the result of Matsuno \cite[Proposition
7.5]{Mat} and Proposition \ref{Selmer general} to obtain similar
results to those in Theorems \ref{Selmer vanishing} and \ref{twist
Selmer vanishing}. Similarly, one can obtain analogous results to
those in Section \ref{Integrality of Akashi series}.

We now state the following analogue of Hachimori's formula
\cite[Theorem 5.4]{HS} (see also \cite[Theorem 16]{Bh},
\cite[Corollary 6.10]{CH}, \cite[Theorem 2.8]{Ho} and \cite[Theorem
3.1]{HV}) for a general Galois representation, and we give a proof
of this statement for completeness.

\bp \label{Selmer general} Assume that $(i)$ $X(A/F^{\cyc})$ is
finitely generated over $\Op$, $(ii)$ $A(F^{\cyc})$ is finite, and
$(iii)$ $F_{\infty}$ is an $S$-admissible $p$-adic extension of $F$
with $G=\Gal(F_{\infty}/F)$ being a pro-$p$ group with no
$p$-torsion. Then $H^2(G_S(F_{\infty}), A)=0$, $\la_S(F_{\infty})$
is surjective and we have the following formula
  \[ \rank_{\Op\ps{H}}X(A/F_{\infty}) = \rank_{\Op}X(A/F^{\cyc}) +
 \sum_w\big(\rank_{\Op}B_w(F^{\cyc}_v)-
 \rank_{\Op\ps{H_w}}B_w(F_{\infty})\big), \]
where $w$ runs through the primes of $F^{\cyc}$ above $S$, $H_w$ is
the decomposition group of $H$ at some prime of $F_{\infty}$ above
$w$, and $B_w$ denotes $A$ or $D_w$ according as $w$ does not or
does divide $p$. In particular, we have
$$\rank_{\Op\ps{H}}X(A/F_{\infty}) \geq \rank_{\Op}X(A/F^{\cyc}).$$
\ep

\bpf Since $X(A/F^{\cyc})$ is finitely generated over $\Op$ and $G$
is pro-$p$ with no $p$-torsion, it follows that $X(A/F^{\cyc})$ is
finitely generated over $\Op\ps{H}$. In particular, we have that
$H^2(G_S(F_{\infty}), A)=0$ and $\la_S(F_{\infty})$ is surjective.
Then we have the following commutative diagram
\[  \entrymodifiers={!! <0pt, .8ex>+} \SelectTips{eu}{}\xymatrix{
    0 \ar[r]^{} & S(A/F^{\cyc}) \ar[d]_{\al} \ar[r] &  H^1(G_S(F^{\cyc}), A)
    \ar[d]_{\be}
    \ar[r] & \bigoplus_{v\in S}J_v(A/F^{\cyc}) \ar[d]_{\ga} \ar[r] & 0 &\\
    0 \ar[r]^{} & S(A/F_{\infty})^H \ar[r]^{} & H^1(G_S(F_{\infty}), A)^H \ar[r] & \
    \bigoplus_{v\in S}J_v(A/F_{\infty})^H \ar[r] &  H^1\big(H, S(A/F_{\infty})\big) \ar[r] & \cdots } \]
with exact rows.  To simplify notation, we write $W_{\infty}=
H^1(G_S(F_{\infty}), A)$ and
 $J_{\infty} = \bigoplus_{v\in S}J_v(A/F_{\infty})$. which in turns yield a long exact
 sequence
 \[ \ba{c} 0\lra \ker\al \lra \ker \be
 \lra \ker \ga
 \lra \coker \al \lra \coker \be \\
  \lra \coker\ga \lra H^1\big(H, S(A/F_{\infty})\big)
 \lra H^1(H, W_{\infty})
 \lra H^1(H, J_{\infty})\lra \cdots \\
 \cdots\lra H^{i-1}(H, J_{\infty}) \lra H^i\big(H, S(A/F_{\infty})\big)
 \lra H^i(H, W_{\infty})
 \lra H^i(H, J_{\infty})\lra \cdots \ea \]
 of cofinitely generated $\Op$-modules. By assumption (i), the
 following exact sequence
\[ 0\lra \ker\al \lra S(A/F^{\cyc})
 \lra S(A/F_{\infty})^H \lra \coker \al \lra 0\] is an exact sequence of cofinitely generated $\Op$-modules.
Comparing the $\Op$-ranks of the terms in both exact sequences, we
obtain
\[ \sum_{i\geq 0}(-1)^i\rank_{\Op}H_i(H,X(A/F_{\infty})) =
\rank_{\Op}X(A/F^{\cyc}) +\sum_{j\geq 1}(-1)^{j}H_j(H,
A(F_{\infty})) \]

\[\hspace{1in} + \sum_w\sum_{j\geq 1}(-1)^{j+1}H_j(H_w,
B_w(F_{\infty})).
\] The required formula then follows from an application of
a formula of Howson \cite[Theorem 1.1]{Ho} on each of the
alternating sum. \epf

\medskip
Another interesting example one may consider is the Galois
representation coming from a primitive Hecke eigenform of weight $k>
2$ for $GL_2 /\Q$, which is ordinary at $p$. Here the finiteness of
$A(F^{\cyc})$ is established in \cite[Proof of Lemma 2.2]{Su}. Also,
as in the case of an elliptic curve, by combining results of Kato
and Hachimori-Ochiai, we can find examples where one can construct
Akashi series for these Selmer groups unconditionally. The
nonexistence of nontrivial pseudo-null submodules for the dual
Selmer groups is established in certain cases \cite[Theorem
3.1]{Su}. As in the case of abelian varieties, one can therefore
obtain results for the Akashi series of these Selmer groups similar
to Theorem \ref{Selmer vanishing}.

\begin{ack}
     This work was written up when the author is a Postdoctoral fellow at the GANITA Lab
    at the University of Toronto. He would like to acknowledge the
    hospitality and conducive working conditions provided by the GANITA
    Lab and the University of Toronto.
        \end{ack}

\end{document}